\documentclass[aos,MSNbibl,nameyear,seceqn,dvips]{arximspdf}
\usepackage{dcolumn}
\usepackage{mathrsfs}
\usepackage{graphicx}

\doi{10.1214/12-AOS1002} 
\volume{40}
\issue{3}
\pubyear{2012}
\firstpage{1285}
\lastpage{1314}

\makeatletter
\newcolumntype{d}[1]{D{.}{.}{#1}}
\newproclaim{Assumption}{Assumption}
\newtheorem{Proposition}{Proposition}
\newtheorem{Theorem}{Theorem}
\newtheorem{Lemma}{Lemma}
\makeatother

\begin{document}
\begin{frontmatter}

\title{Test for bandedness of high-dimensional covariance matrices and bandwidth estimation}
\runtitle{Bandedness test for covariance matrices}

\begin{aug}
\author[A]{\fnms{Yumou} \snm{Qiu}\ead[label=e2]{yumouqiu@iastate.edu}}
\and
\author[B]{\fnms{Song Xi} \snm{Chen}\corref{}\thanksref{aut1}\ead[label=e1]{songchen@iastate.edu}}
\runauthor{Y. Qiu and S. X. Chen}
\thankstext{aut1}{Supported by NSFC key Grant 11131002.}
\affiliation{Iowa State University, and Peking University and Iowa State University}
\address[B]{Department of Statistics\\
Iowa State University\\
Ames, Iowa 50011-1210\\
USA\\
\printead{e2}}
\address[A]{Department of Business Statistics and Econometrics\\
Guanghua School of Management and\\
\quad Center for Statistical Science\\
Peking University\\
Beijing 100871\\
China\\
and\\
Department of Statistics\\
Iowa State University\\
Ames, Iowa 50011-1210\\
USA\\
\printead{e1}}
\end{aug}

\received{\smonth{11} \syear{2011}}
\revised{\smonth{4} \syear{2012}}

%
\begin{abstract}
Motivated by the latest effort to employ banded matrices to estimate a
high-dimensional
covariance $\Sigma$, we propose a test for $\Sigma$ being banded with
possible diverging bandwidth. The test is adaptive to the ``large $p$,
small $n$'' situations
without assuming a specific parametric distribution for the data. We
also formulate a consistent estimator for the bandwidth of a banded
high-dimensional covariance matrix. The properties of the test and the
bandwidth estimator are investigated by theoretical evaluations and
simulation studies, as well as an empirical analysis on a protein mass
spectroscopy data. 
\end{abstract}

%

\end{frontmatter}

\section{Introduction}\label{sec1}

High-dimensional data are increasingly collected in statistical
applications, which include biological experiments, climate and
environmental studies, financial observations and others. The high
dimensionality calls for new statistical methodologies which are
adaptive to this new feature of the modern statistical data.
The covariance matrix $\Sigma=\operatorname{Var}(X)$ for a
$p$-dimensional random
vector $X$ is an important measure on the dependence among components
of $X$.
The sample covariance $S_n$, constructed based on $n$ independent
copies of $X$, 
is a key ingredient in many statistical procedures in the conventional
multivariate analysis [\citet{r2} and \citet{r24}] where the
data dimension $p$ is regarded as fixed.
The widespread use of $S_n$ in the conventional multivariate procedures
is largely due to $S_n$ being a consistent estimator of $\Sigma$ when
$p$ is fixed or small relative to the sample size $n$.
However, for high-dimensional data such that $p/n \to c \in(0, \infty
]$, it is known that the eigenvalues of the sample covariance matrix
are no longer consistent to their population counterpart, as
demonstrated in Bai and Yin (\citeyear{r5}), Bai, Silverstein and Yin (\citeyear{r4}),
Johnstone (\citeyear{J01}) and El~Karoui (\citeyear{r15}). These mean that the sample
covariance $S_n$ is no longer consistent to $\Sigma$, which hinders
applications of many conventional multivariate statistical procedures
for high-dimensional data.

To overcome the problem with the sample covariance, constructing
covariance estimators via banding or tapering
the sample covariance matrix has been a focus in high-dimensional
covariance estimation. 
\citet{r31} considered banding the Cholesky factor matrix via the
kernel smoothing estimation, {which was further developed by
\citet{r27}}. 
\citet{r7} proposed banding the sample covariance matrix directly
for estimating $\Sigma$ and banding the Cholesky factor matrix for
estimating $\Sigma^{-1}$. They demonstrated that both estimators are
consistent to~$\Sigma$ and $\Sigma^{-1}$, respectively, for some
``bandable'' classes of covariance matrices. \citet{r12} proposed a
tapering estimator, which can be viewed as a soft banding on the sample
covariance, which was designed to improve the banding estimator of
Bickel and Levina. They demonstrated that the tapering estimator
attains the optimal minimax rates of convergence for estimating the
covariance matrix. Wagaman and Levina (\citeyear{r30}) developed a method for
discovering meaningful orderings of variables such that banding and
tapering can be applied. Both the banding and tapering methods for
covariance estimation are well connected to the regularization method
considered in Huang et al. (\citeyear{r18}), \citet{r8},
\citet{r16} and \citet{r26}.

Motivated by the promising results regarding banding and tapering the
sample covariance, we develop in this paper a test procedure on the
hypothesis that $\Sigma$ is banded. The rationale for developing such a
test is to check a $\Sigma$ in the so-called ``bandable'' class outlined
in \citet{r7} such that the banding or the tapering estimators are
consistent. 
There is yet a practical guideline to confirm or otherwise if a $\Sigma
$ is within the ``bandable'' class so that the banding and tapering can
be applied.
Hence, a direct testing on $\Sigma$ being banded provides a path of
advance to gain knowledge on the structure of the covariance. If the
banded hypothesis is confirmed by the test, the banding and tapering
estimators may be employed.

Diagonal matrices are the simplest among banded matrices.
Given the importance commanded by covariance matrices in
high-dimensional multivariate analysis, directly testing for $\Sigma$
being diagonal and the so-called sphericity hypothesis in classical
multivariate analysis [\citet{r20} and \citet{r25}], have
been considered in a set of studies including \citet{r21},
\citet{r19}, \citet{r28}, \citet{r13}
and Cai and Jiang (\citeyear{r10}) under high dimensionality. For normally
distributed data, \citet{r19} proposed testing for diagonal
$\Sigma
$ by considering a coherence statistic $L_{n}=\mathop{\max }_{1\leq i<j\leq
p}|\hat{\rho}_{ij}|$, where $\hat{\rho}_{ij}$ is
the sample correlation coefficient between the $i$th and the $j$th
components of the random vector $X$. Jiang established the asymptotic
distribution of $L_{n}$ under the null diagonal hypothesis, which was
used to derive a sphericity test. As $L_n$ is an extreme value type,
its convergence to its limiting distribution can be slow.
\citet{r23} proposed a modification which is shown to be able to
speed up the convergence.
Cai and Jiang (\citeyear{r10}) extended the test of \citet{r19} for the
bandedness of $\Sigma$, which is shown to be applicable for the ``large
$p$, small $n$'' situations such that $\log(p)=o(n^{1/3})$.

In this paper, we propose a nonparametric test for $\Sigma$ being
banded without assuming a parametric distribution for the
high-dimensional data. The test is formulated to allow the dimension to
be much larger than the sample size. 
Based on the test statistic for bandedness, we propose a consistent
estimator for the bandwidth of a banded high-dimensional covariance.
The properties of the test and bandwidth estimator are demonstrated by
theoretical evaluation, simulation studies and empirical analysis on a
protein mass spectroscopy data for prostate cancer.

The paper is organized as follows. Section \ref{sec2} introduces the hypotheses,
the assumptions and the test statistic. In Section \ref{sec3}, we present the
properties of the test statistic and the test, and evaluate its power
properties. Estimation of the bandwidth is considered in Section \ref{sec4}.
Section \ref{sec5} reports simulation results. An empirical analysis on a
prostate cancer spectroscopy data is outlined in Section \ref{sec6}. All
technical details are relegated to the
\hyperref[app]{Appendix}.\looseness=1

\section{Preliminary}\label{sec2}

Let $X_{1},X_{2},\ldots, X_{n}$ be independent and identically
distributed $p$-dimensional random vectors with mean $\mu$ and
covariance matrix $\Sigma=(\sigma_{ij})_{p\times p}$. A matrix
$A=(a_{ij})_{p \times p}$ is said to be banded if there exists an
integer $k \in\{ 0,\ldots, p-1\}$
such that $a_{ij}=0$ for $|i-j|>k$. The smallest~$k$ such that $A$ is
banded is called the bandwidth of $A$. Banding of $A$ at a~bandwidth~$k$
refers to setting $a_{ij} =0$
for all $|i-j| > k$. 

Let $B_{k}(\Sigma)=(\sigma_{ij}\mathrm{I}\{|i-j|\leq k\})_{p\times p}$ be
a banded version of $\Sigma$ with bandwidth~$k$. Specifically,
$B_{0}(\Sigma)$ is the diagonal version of $\Sigma$.
We intend to test
%
\begin{equation} \label{eq:hypo}
H_{k,0}\dvtx \Sigma=B_{k}(\Sigma)\quad\mbox{vs.}\quad
H_{k,1}\dvtx\Sigma\neq B_{k}(\Sigma)
\end{equation}
for $k=o(p^{1/4})$. 
Hence, the bandwidth $k$ of $\Sigma$ to be tested can be either fixed
or diverging to infinite as long as it is slower than $p^{1/4}$.
Allowing divergent bandwidth in the hypothesis is an improvement over
the sphericity test 
as considered in \citet{r21} and \citet{r13}.
It also connects to the latest works on high-dimensional covariance
estimation with banded or tapered versions of the sample covariance as
in \citet{r7} and \citet{r12}. In particular, \citet{r12} showed that the optimal minimax rates for the bandwidth of the
banded covariance estimator of \citet{r7} is $k=O[ \{n/\log(p)\}
^{1/(2 \alpha+1)}]$, and that for the tapering estimator is
$k=O(n^{1/(2 \alpha+1)})$,
where $\alpha$ is an index value for a ``bandable''\vadjust{\goodbreak} class of covariances
\begin{eqnarray}\label{eq:bandable}
\mathfrak{U}(\varepsilon_{0},\alpha,C)&=& \biggl\{\Sigma\dvtx\max
_{j}\sum_{|i-j|>k} |\sigma_{ij}| \leq C k^{-\alpha}\mbox{ for all }k>0,\nonumber\\[-8pt]\\[-8pt]
&&\hspace*{12pt}\mbox{and  }0<\varepsilon_{0}\leq\lambda_{\min}(\Sigma)\leq
\lambda
_{\max}(\Sigma)\leq\varepsilon_{0}^{-1}  \biggr\}.\nonumber
\end{eqnarray}
%
The range of bandwidths $k=o(p^{1/4})$ in the hypothesis (\ref{eq:hypo}) should cover the above optimal rates when $p \gg  n$.

We note that $H_{k,0}$ is valid if and only if $\sum_{ |i-j| > k_{p}}
\sigma_{ij}^2 =0$, and the latter implies that $\operatorname{tr}\{
\Sigma
-B_{k}(\Sigma)\}^2=0$. A strategy is to construct an unbiased estimator
of $\operatorname{tr}\{\Sigma-B_{k}(\Sigma)\}^2$ and use it to
develop the test
statistic. Let $D_{q}:=\sum_{l=1}^{p-q}\sigma_{l l+q}^{2}$ be the sum
of squares of the $q$th sub-diagonal of $\Sigma$.
Then, $\operatorname{tr}\{\Sigma-B_{k}(\Sigma)\}^2 = 2\sum_{q
=k+1}^{p-1} D_{q}$.
It can be checked that an unbiased estimator of $D_{q}$ is 
\begin{eqnarray*}
\hat{D}_{n q}&=&\sum _{l=1}^{p-q} \Biggl\{\frac
{1}{P^{2}_{n}}\sum
 _{i,j}^{*}(X_{il}X_{i l+q})(X_{jl}X_{j l+q})
-2\frac{1}{P^{3}_{n}}\sum
_{i,j,k}^{*}X_{il}X_{k  l+q}(X_{jl}X_{j l+q})\\
&&\hspace*{158pt}{}+\frac{1}{P^{4}_{n}}\sum
_{i,j,k,m}^{*}X_{il}X_{j l+q}X_{kl}X_{m l+q} \Biggr\},
\end{eqnarray*}
where $\sum ^{*}$ denotes summation over mutually different
subscripts shown
and $P_{n}^{b} = n ! /(n - b) !$. The reason to sum over different
indices is for easier manipulations with the mean and variance of the
final test statistic and to establish the asymptotic normality. The
latter leads to a test procedure for the bandedness.

We consider the following statistic:
%
\begin{equation}\label{eq:Wnk0}
W_{nk}:= 2\sum_{q=k+1}^{p-1} \hat{D}_{n q}.
\end{equation}
%
As each $\hat{D}_{n q}$ is invariant under the location shift, $W_{nk}$
is also location shift invariant. Hence, without loss of generality, we
assume $\mu={\mathrm{E}}(X)=0$.

To facilitate our analysis, as \citet{r3} and \citet{r13},
we assume a multivariate model for the high-dimensional data.

\begin{Assumption}\label{ass1}
(i) $X_1, X_2, \ldots, X_n$ are independent and identically distributed (i.i.d.)
$p$-dimensional random vectors such that
%
\begin{equation}\label{eq:ass1a}
X_i=\Gamma Z_i\qquad \mbox{for $i=1,2,\ldots,n$,}
\end{equation}
where
$\Gamma$ is a $p\times m$ constant matrix with
$m\ge p$, $\Gamma\Gamma' = \Sigma$, and $Z_1,\ldots,Z_n$ are i.i.d. $m$-dimensional random vectors such that
${\mathrm{E}}(Z_1)=0$ and
${\operatorname{Var}}(Z_1)=I_{m}$.\vadjust{\goodbreak}

(ii) Write $Z_1=(z_{1 1}, \ldots,
z_{1 m})^T$. Each $z_{1 l}$ has uniformly bounded
$8$th moment, and there exist finite constants $\Delta$ and $\omega$ such
that for $l=1,
\ldots, m$, ${\mathrm{E}}(z_{1 l}^4)=3+\Delta$, ${\mathrm{E}}(z_{1
l}^3)=\omega$ %
and for any integers $\ell_\nu\ge0$ with
$\sum_{\nu=1}^q\ell_\nu=8$
%
\begin{equation}\label{eq:mixedmoment}
{\mathrm{E}}(z_{i_1}^{\ell_1}z_{i_2}^{\ell_2}\cdots
z_{i_q}^{\ell_q})={\mathrm{E}}(z_{1 i_1}^{\ell_1}){\mathrm{E}}(z_{1
i_2}^{\ell_2})
\cdots{\mathrm{E}}(z_{1 i_q}^{\ell_q})
\end{equation}
whenever
$i_1,i_2,\ldots,i_q$ are distinct subscripts.
\end{Assumption}

The requirement of common third and fourth moments of $z_{1l}$ is not
essential and is purely for the sake of simpler notation. Our theory
allows different third and fourth moments as long as they are uniformly
bounded, which are actually assured by $z_{1l}$ having uniformly
bounded $8$th moment.

The asymptotic framework that regulates the sample size $n$, the
dimensionality $p$ and the covariance $\Sigma$ is the following.

\begin{Assumption}\label{ass2}
As $n\to\infty$, $p=p(n)\to\infty$,
$n =O(p)$ 
and ${\operatorname{tr}(\Sigma^4)}/\break{\operatorname{tr}^2(\Sigma
^2)}=O({p}^{-1})$.
\end{Assumption}

We note that $n=O(p)$ includes $p \gg  n$, the ``large $p$, small $n$''
paradigm, but may not imply $p=O(n)$. Different from the usual approach
of specifying an explicit growth rate of $p$ with respect to $n$,
Assumption \ref{ass2} requires ratio of $\operatorname{tr}(\Sigma^4)$ to
${\operatorname{tr}
}^{2}(\Sigma
^2)$ shrinks at the rate of $p^{-1}$ or smaller.
The latter is stronger than ${\operatorname{tr}(\Sigma
^4)}/{\operatorname{tr}^2(\Sigma
^2)}=o(1)$. It is needed due to possible diverging bandwidths.

Let
\[
\mathcal{U}_p= \biggl\{\Sigma\dvtx \frac{\operatorname{tr}(\Sigma
^4)}{\operatorname{tr}^2(\Sigma
^2)}=O(p^{-1}) \biggr\}
\]
be the class of covariances satisfying the last part of Assumption \ref{ass2}.
The class includes the ``bandable'' class $\mathfrak{U}(\varepsilon
_{0},\alpha,C)$ of \citet{r7} given in (\ref{eq:bandable})
for the banding estimation.
To appreciate this, let $\lambda_1\leq\lambda_2\leq\cdots\leq
\lambda
_p$ be the eigenvalues of $\Sigma$. If the smallest and largest
eigenvalues are bounded away from 0 and $\infty$ respectively, then
\[
\frac{\operatorname{tr}(\Sigma^4)}{\operatorname{tr}^{2}(\Sigma
^2)}=\frac{\sum
_{i=1}^{p}\lambda
_{i}^{4}}{ (\sum_{i=1}^{p}\lambda_{i}^{2} )^{2}}
\leq\frac{\lambda_{p}^{4}}{p\lambda_{1}^{4}}=O(p^{-1}).
\]
Therefore, the ``bandable'' covariances are contained in $\mathcal
{U}_p$. Now suppose that $\Sigma$ has exactly $m_p$ zero eigenvalues
and $\lambda_{m_{p}+1}$ being the smallest nonzero eigenvalue. 
Then
\[
\frac{\operatorname{tr}(\Sigma^4)}{\operatorname{tr}^{2}(\Sigma
^2)} 
\leq\frac{\lambda_{p}^{4}}{(p-m_{p})\lambda_{m_{p}+1}^{4}}.
\]
Thus, $\Sigma$ is in $\mathcal{U}_p$ as long as
$\lambda_{p}/\lambda_{m_p +1}$ is bounded and $m_{p}\leq cp$ for some
\mbox{$c\in(0,1)$} as $p \to\infty$. The latter means that the class
$\mathcal
{U}_p$ is likely to contain the class considered in \citet{r12},
which allows the smallest eigenvalue to diminish to zero.
It can be also checked that the following two covariances,
\[
\Sigma= \bigl(\sigma_i \sigma_j \rho^{|j-i|} \bigr)_{p\times p}\quad\mbox{or}\quad
\Sigma= \bigl(\sigma_i \sigma_j \rho^{|j-i|} \mathrm{I}( |j-i| \le d)
\bigr)_{p\times p},
\]
are members of $\mathcal{U}_p$ if
$\{\sigma_l^2\}_{l=1}^p$ are uniformly bounded from infinity and zero
respectively.

\section{Main results}\label{sec3}

We first describe the basic properties of the statistic~$W_{n k}$
defined in (\ref{eq:Wnk0}). Let
\begin{eqnarray}\label{eq:nu}
\nu_{nk}^{2}&=&\frac{4}{n^2}\operatorname{tr}^2(\Sigma^2)+\frac
{8}{n}\operatorname{tr} \bigl\{
\Sigma\bigl(\Sigma-B_{k}(\Sigma)\bigr) \bigr\}^2\nonumber\\[-8pt]\\[-8pt]
&&{}+\frac{4}{n}\Delta\operatorname{tr}\bigl\{\Gamma'\bigl(\Sigma
-B_{k}(\Sigma)\bigr)\Gamma
\circ
\Gamma'\bigl(\Sigma-B_{k}(\Sigma)\bigr)\Gamma\bigr\},\nonumber
\end{eqnarray}
where $\Omega\circ\Lambda=(\omega_{ij}\lambda_{ij})$ for two matrices
$\Omega=(\omega_{ij})$ and $\Lambda=(\lambda_{ij})$.

\begin{Proposition}\label{pro1}
 Under Assumptions \ref{ass1} and \ref{ass2},
\[
{\mathrm{E}}(W_{n k})=\operatorname{tr}[\{\Sigma- B_{k}(\Sigma)\}
^2] \quad\mbox{and}\quad
 {\operatorname{Var}}(W_{nk})=\nu_{nk}^{2}+o(\nu_{nk}^{2})
.
\]
\end{Proposition}

  The proposition indicates that under $H_{k,0}$,
\[
   {\mathrm
{E}}(W_{n k})=0\quad
\mbox{and} \quad\nu_{nk}=2\operatorname{tr}[\{B_{k}(\Sigma)\}^{2}]/{n},
\]
 and
$\nu_{n k}^{2}$ is the leading order variance of $W_{nk}$. It can be
shown that
$\operatorname{tr}\{\Sigma(\Sigma-B_{k}(\Sigma))\}^{2}\leq
4(k+1)^{2}{\operatorname{tr}
}(\Sigma
^4)$. Since
\[
\operatorname{tr}\bigl\{\Gamma'\bigl(\Sigma-B_{k}(\Sigma)\bigr)\Gamma\circ
\Gamma'\bigl(\Sigma
-B_{k}(\Sigma)\bigr)\Gamma\bigr\}
\leq\operatorname{tr}\bigl\{\Sigma\bigl(\Sigma-B_{k}(\Sigma)\bigr)\bigr\}^{2},
\]
$\Delta\geq-2$ and ${\operatorname{tr}(\Sigma
^{4})}/{\operatorname{tr}^{2}(\Sigma
^{2})}=O(p^{-1})$, we have
%
\begin{equation}\label{eq:bounds}
4n^{-2}\operatorname{tr}^2(\Sigma^2)\leq\nu_{nk}^{2}\leq
C_{0}a_{np}{\operatorname{tr}
}^2(\Sigma
^2)
\end{equation}
for a constant $C_{0}\geq4$ and $a_{np}=n^{-2}+k^{2}(np)^{-1}$. We note
that $a_{np}\to0$ as $n\to\infty$ since $k=o(p^{1/4})$. In particular,
if $k$ is fixed, $a_{np}=O(n^{-2})$.

The following theorem establishes the asymptotic normality of $W_{nk}$.

\begin{Theorem}\label{teo1}
Under Assumptions \ref{ass1} and \ref{ass2}, and if $k=o(p^{1/4})$,
\[
\frac{W_{nk}-\operatorname{tr}[\{\Sigma-B_{k}(\Sigma)\}^{2}]}{\nu
_{nk}}\stackrel{D}{\rightarrow}
N(0,1).
\]
\end{Theorem}

  In order to formulate a test procedure based on the asymptotic
normality, we need to estimate $\operatorname{tr}[\{B_{k}(\Sigma)\}
^{2}]$ since
$\nu
_{nk}=2\operatorname{tr}[\{B_{k}(\Sigma)\}^{2}]/{n}$ under~$H_{k, 0}$.
Let $V_{nk}:= \hat{D}_{n 0}+2\sum_{q=1}^{k} \hat{D}_{n q}$ be the
estimator, whose consistency to $\operatorname{tr}[\{B_{k}(\Sigma)\}
^{2}]$ is
implied in the following proposition.\vadjust{\goodbreak} 

\begin{Proposition}\label{pro2}
Under Assumptions \ref{ass1} and \ref{ass2}, ${\operatorname
{Var}}\{
{V_{nk}}/{\operatorname{tr}(\Sigma^{2})}\}=O(a_{np})$, where
$a_{np}=n^{-2} +
k^{2}(np)^{-1}$.
\end{Proposition}

Since ${\mathrm{E}}(V_{nk})=\operatorname{tr}[\{B_{k}(\Sigma)\}
^{2}]$ and
$a_{np}\to0$, Proposition \ref{pro2} means that, under $H_{k,0}$,
${V_{nk}}/{\operatorname{tr}[\{B_{k}(\Sigma)\}^{2}]}\stackrel{p}{\rightarrow}1$ as $n\to
\infty$.
This together with Theorem \ref{teo1}
indicates that under $H_{k,0}$
\[
T_{n k} =: n\frac{W_{n k}}{V_{nk}}\stackrel{D}{\rightarrow}N(0,4).
\]
This leads to our choice of $T_{nk}$ as the test statistic and the
proposed test of size $\alpha$ that rejects $H_{k,0}$ if $T_{n k} \ge2
z_{\alpha}$ where $z_{\alpha}$ is the upper $\alpha$ quantile of $N(0,1)$.

As Theorem \ref{teo1} prescribes the asymptotic normality under both
$H_{k,0}$\break
and~$H_{k,1}$, it permits a power evaluation of the test.
Let
%
\begin{equation}
\delta_{nk}=\frac{\operatorname{tr}(\Sigma^2)-\operatorname
{tr}[\{B_{k}(\Sigma)\}
^{2}]}{\nu_{nk}},
\end{equation}
which may be viewed as a signal to noise ratio for the testing problem.
This is because $\operatorname{tr}[\{ \Sigma- B_{k}(\Sigma)\}
^{2}]$ is the square
of Frobenius norm of the difference between $\Sigma$ and its $k$-banded
version, and ${\nu_{nk}}$ measures the level of noise in the statistic~$W_{n k}$.
Then, the power of the test under $H_{k,1}\dvtx \Sigma\ne B_{k}(\Sigma)$ is
\begin{eqnarray*}
\beta_{nk} &=& P\{n{W_{nk}}/{V_{nk}}\ge2z_{\alpha} | \Sigma\ne
B_{k}(\Sigma)\} \\
&=& P \biggl(\frac{W_{nk}-\operatorname{tr}(\Sigma^2)+{\operatorname
{tr}}[\{B_{k}(\Sigma)\}
^{2}]}{\nu
_{nk}}\geq\frac{2 z_{\alpha}V_{nk}}{n\nu_{nk}}-\delta_{nk} \biggr).
\end{eqnarray*}
Since $\nu_{nk}\geq2n^{-1}\operatorname{tr}(\Sigma^{2})$, then
$2V_{nk}/(n\nu
_{nk})\leq V_{nk}/\operatorname{tr}(\Sigma^{2})$ for $n$ large. Hence asymptotically,
%
\begin{equation}\label{eq:beta}
\beta_{nk}\geq P
\biggl(\frac{W_{nk}-\operatorname{tr}(\Sigma^2)+\operatorname{tr}[\{
B_{k}(\Sigma)\}^{2}]}{\nu
_{nk}}\geq z_{\alpha}\frac{V_{nk}}{\operatorname{tr}(\Sigma
^{2})}-\delta_{nk}
\biggr).
\end{equation}

To gain more insight on the power, let $r_{k}=\operatorname{tr}[\{
B_{k}(\Sigma)\}
^{2}]/\operatorname{tr}(\Sigma^2)$.
Clearly, $r_{k} \le1 $ and is monotone nondecreasing with respect to
$k$. If $\Sigma$ is banded with bandwidth~$k_0$, then
%
\begin{equation}\label{eq:cpoint}
r_{k} < 1\qquad  \mbox{for }k < k_0\quad  \mbox{and}\quad
  r_{k} =
1 \qquad \mbox{for }k \ge k_0.
\end{equation}
From the bounds for $\nu_{nk}$ in (\ref{eq:bounds}), it follows that
%
\begin{equation}\label{eq:bound2}
(C_{0}a_{np})^{-1/2} (1-r_{k}) \leq\delta_{nk} \leq {\tfrac
{1}{2}}{n}(1-r_{k}),
\end{equation}
which indicates that $a_{np}^{-1/2}(1-r_{k})=O(\delta_{nk})$. When
$k$ is fixed, $a_{np}=O(n^{-2})$ and $\delta_{nk}\sim n(1-r_{k})$,
indicating that $\delta_{n k}$ is at the exact order of $n (1-r_k)$.


\begin{Theorem}\label{teo2}
Under Assumptions \ref{ass1} and \ref{ass2}, $H_{k,1}$ and if
$k=o(p^{1/4})$, then:
\begin{longlist}
\item $\liminf _n\beta_{nk} \ge1-\Phi(z_{\alpha}-\liminf_n\delta_{nk} )$;

\item if $a_{np}^{-1/2}(1-r_{k})\to\infty$, then $\beta_{nk}\to1$ as
$n\to\infty$.\vadjust{\goodbreak}
\end{longlist}
\end{Theorem}

Theorem \ref{teo2} indicates that the proposed test is consistent as long as the
speed of $1-r_{k} \to0$ under $H_{k,1}$ is not faster than
$a_{np}^{1/2}$. The test will have nontrivial power as long as
$\liminf
 _n\delta_{nk} > 0$. 
If $n (1-r_{k}) \to0$, the test will have no power beyond the
significant level $\alpha$. We\vspace*{1pt} note that this happens when $H_{k,0}$
and $H_{k,1}$ are extremely close to each other, so that $1-r_k$ decays
to zero faster than
$n^{-1}$. We are actually a little amazed by the fact that the test is
powerful as long as \mbox{$\liminf _n a_{n p}^{-1/2} (1- r_{k}) > 0$}
or equivalently
$(1- r_{k})$ does not shrink to zero faster than $a_{n p}^{1/2}$,
despite the high dimensionality and a possible diverging bandwidth $k$.
Theorem \ref{teo2} and (\ref{eq:bound2}) together imply that if $r_{k}$ does not
vary much as $p$ increases, the power of the test will be largely
determined by $n$, as confirmed by our simulation study in Section \ref{sec5}.

Our proposed test is targeted on the covariance matrix $\Sigma$. A test
for the correlation matrix can be developed by modifying the test
statistic by first standardizing each data dimension via its sample
standard deviation.
The theoretical justification would be quite involved, and would
require extra effort.
In addition to be invariant under the location shift, the test
statistic is invariant if all the variables among the high-dimensional
data vector are transformed by a common scale. However, the proposed
test statistic is not invariant under variable-specific scale
transformation. The above mentioned test for the correlation matrix
would be invariant under variable-specific scale transformation.

\section{Bandwidth estimation}\label{sec4}

We propose in this section an estimator to the bandwidth of banded
covariance $\Sigma$.
Estimating the bandwidth of a banded covariance matrix is an important
and practical issue, given the latest advances on covariance estimation
by banding [\citet{r7}] or tapering [\citet{r12}] sample
covariance matrices. Indeed, finding an adequate bandwidth is a
pre-requisite for applying either the banding or tapering estimators.

The proposed estimator is motivated by the test procedure developed in
the previous section.
Let $k_0$ be the true bandwidth.
As the proposed test is consistent as long as $r_{k} \to1$ not too
fast, and the sample size is large enough (can still be much less than
$p$), the proposed test would reject (not reject) $H_{k, 0}$ for $k$
less (larger) than~$k_0$.
An immediate but rather naive strategy would be to use the smallest
integer $k$ such that $H_{k,0}$ is not rejected as the bandwidth estimator.
However, this strategy may be insufficient to counter ``abnormal''
samples which can produce larger (smaller) values of the statistic
$\tilde{T}_{nk}:=W_{nk}/V_{n k}$ consistently for a wide range of $k$
values, when in fact $H_{k,0}$ ($H_{k,1}$) is true. And yet these
``abnormal'' samples are expected within the normal range of variations.
To make the estimator robust against these ``abnormal'' samples and not
so much dependent on the significant level $\alpha$, we consider an
estimator based on the difference between successive statistics,
$d_{nk}=\tilde{T}_{nk}-\tilde{T}_{n k+1}$.\vadjust{\goodbreak}


We assume the true bandwidth $k_0$ be either fixed or diverging as long as
%
\begin{equation}\label{eq:k0}
k_0 (n^{-1/2} + p^{-1/4}) \to0\qquad  \mbox{as }n \to
\infty,
\end{equation}
which covers a quite wide range for the bandwidth.
Note that
%
\[
\tilde{T}_{nk} = \frac{W_{nk}-{\mathrm{E}}(W_{nk})}{\nu_{nk}}
\frac{ \nu_{nk}}{V_{nk}}+ \frac{ {\mathrm{E}}(W_{n k}) }{V_{nk}}.
\]
For $k\leq M$, where $M =o(p^{1/4})$ is a pre-chosen sufficiently large
integer, ${\{W_{nk}-{\mathrm{E}}(W_{nk})\}}/{\nu_{nk}}$ is stochastically
bounded (Theorem \ref{teo1}) and
from (\ref{eq:bounds}), we have
\[
\tilde{T}_{nk} = O_p \biggl( \frac{a_{n p}^{{1}/{2}} r_{k}^{-1}
\operatorname{tr}[\{B_{k}(\Sigma)\}^{2}] }{V_{n k} } \biggr) +
(r_{k}^{-1} - 1)
\frac{ \operatorname{tr}[\{B_{k}(\Sigma)\}^{2}]}{V_{n k} }.
\]
Let $b_{n k}=V_{n k}/\operatorname{tr}[\{B_{k}(\Sigma)\}^2]-1$.
From Propositions \ref{pro1}
and \ref{pro2},
%
\begin{equation}
{\mathrm{E}}(b_{n k})=0\quad\mbox{and}\quad
{\operatorname{Var}}(b_{n k})=O(a_{np}r_{k}^{-2}).
\end{equation}
Since $\Sigma=B_{k_0}(\Sigma)$ is nonnegative definite,
${\operatorname{tr}
}(\Sigma
^{2})\leq(2k_{0}+1)\operatorname{tr}[\{B_{0}(\Sigma)\}^{2}]$.
Hence, for any $k$, $r_{k}\geq(2k_{0}+1)^{-1}$. 
These imply that
%
\begin{eqnarray}\label{eq:Wnk}
\tilde{T}_{nk} & = & O_p( a_{n p}^{ {1}/{2} } k_0) + (r_{k}^{-1} -
1) \{ 1 + o_p(1)\}.
\end{eqnarray}
It can be checked that $a_{n p}^{ {1}/{2} } k_0 \to0$ under (\ref
{eq:k0}), which makes the first term on the right of the above equation
negligible relative to the second term. And the second term is quite
indicative between $k < k_0$ and $k \ge k_0$, since $r_{k}=1$ for $k
\ge k_0$.

To amplify the second term when $k < k_0$ while not inflicting the
first term on the right of (\ref{eq:Wnk}) too much, we consider
multiplying $n^{\delta}$ on $\tilde{T}_{nk}$ for a small positive
$\delta$ and let
$ d_{nk}^{(\delta)} = n^{\delta} ( \tilde{T}_{nk} - \tilde{T}_{n k+1})$.
The proposed bandwidth estimator is
%
\begin{equation}\label{eq:fixed}
\hat{k}_{\delta,\theta}=\min\bigl\{k\dvtx \big|d_{n k}^{(\delta)}\big|< \theta\bigr\}
\end{equation}
for a pair of tuning parameters $\delta>0$ and $\theta>0$.
The following theorem gives the consistency of the bandwidth estimator
for both fixed or diverging $k_0$.

%

\begin{Theorem}\label{teo3}
Under Assumptions \ref{ass1} and \ref{ass2}, 
if $\liminf_n \{\inf_{k<k_{0}}(r_{k+1}-r_{k})\}>0$,
then for any $\theta>0$, $\hat{k}_{\delta,\theta} -
k_{0}\stackrel{p}{\rightarrow} 0 $ under either of the two settings:
\textup{(i)} for any $\delta\in(0,1)$ if $k_0$ is bounded; \textup{(ii)} for any
$\delta\in(0,1/2]$ if $k_0$ is diverging but satisfies (\ref{eq:k0}),
and $\{\sigma_{ll}\}_{l=1}^{p}$ are uniformly bounded away from $0$ and
$\infty$.
\end{Theorem}

\begin{figure}

\includegraphics{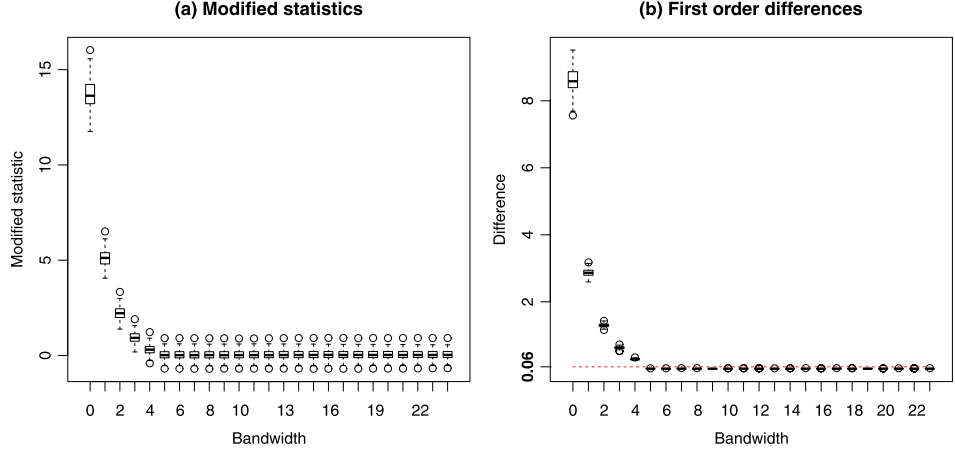}

\caption{Box-plots of the modified statistics $n^{\delta}\tilde
{T}_{nk}$ and their first order differences of the simulated data. The
dashed line in the right panel is $\theta=0.06$. The true bandwidth is 5.}\label{fig1}
\end{figure}

We would like to remark that the multiplier $n^{\delta}$ in
$d_{nk}^{(\delta)}$'s formation leads to~$\theta$ being ``free ranged''
as long as $\theta>0$. If such multiplication is not administrated,
namely by setting $\delta=0$, the range of $\theta$ needs to be
restricted properly to ensure convergence.
The requirement of $\liminf_n \{\inf_{k<k_{0}}(r_{k+1}-r_{k})\}>0$ is
to avoid situations where $\Sigma$ has segments of zero sub-diagonals
followed by nonzero sub-diagonals when one moves away from the main
diagonal. Our estimator can be modified to suit such situations.
However, we would not elaborate here for the sake of simplicity in the
presentation.
Attaining the consistency of $\hat{k}_{\delta, \theta}$ with diverging
$k_0$ requires a smaller $\delta$ value.

To better understand the theorem and the bandwidth estimator, we
conducted a simulation study for $k_0=5$,
$n=60$ and $p=600$ with $X_i$ generated from Model~(\ref{eq:5.1}) with a
multivariate normal distribution. The detailed simulation setting will
be provided in Section \ref{sec5}. Figure \ref{fig1} presents box-plots of the modified
statistics $n^{\delta}\tilde{T}_{nk}$ (left panel) and its first-order
difference $d_{nk}^{(\delta)}$ (right panel), with $\delta=0.5$. We see
from the right panel that the first five boxes are relatively large,
and $d_{nk}^{(\delta)}$ is close to 0 while for $k \ge5$. This
indicates that five would be the bandwidth estimate.

In practical implementations with finite samples, the bandwidth
estimator may be sensitive to the tuning parameters $\delta$ and
$\theta
$. Note that, as revealed a few paragraphs earlier, $d_{nk}$ should be
significantly larger than $0$ for $k<k_{0}$ and close to $0$ for $k\geq
k_{0}$. Such a pattern, as displayed in Figure \ref{fig1}, indicates that $k_0$
is a change point for $\{ d_{nk} \}_{k=0}^{M}$.
This motivates us to consider a~regression change-point detection
algorithm for bandwidth estimation.
Consider $d_{n j}$, the difference between successive statistics $T_{n
j}$, for $j=1,\ldots, M$, for a sufficiently large $M$ that covers
the true bandwidth $k_0$.
The idea is to fit, at each candidate $k$, a regression function
$g_{k}(j)$ to $\{ d_{n j} \}_{j=0}^M$ such that $g(j) \equiv g(k)$ for
all $j > k$.
We may fit a nonparametric, locally weighted linear regression
[\citet{r14}; Fan and Gijbels (\citeyear{r17})] on $j \in L_k =\{l\dvtx 0\leq
l\leq k\}$ to the left of $k$ with the smoothing window-width $hk$,
where $h$ is a smoothing parameter, and fit a flat line at the level
$d_{n k}$ for $j \in R_k =\{l\dvtx k+1\leq l\leq M \}$ to the right of $k$.
If $k$ is too small for the above nonparametric regression, a
parametric polynomial regression may be conducted. Let $\hat{g}_{k}(j)$
be the regression estimate, nonparametric or parametric, obtained over
the set $L_k$, and let
\[
\operatorname{err}(k) = \sum_{j \in L_{k}}|\hat{g}_{k}(j)-d_{nj}|+\sum_{j \in
R_{k}}|d_{n k}-d_{nj}|
\]
be the absolute deviation of the fitted errors. Then a bandwidth
estimator, as we call the change-point estimator, is
%
\begin{equation}\label{eq:empirical}
\hat{k}=\mathop{\operatorname{arg\, min}}_{k} \{\operatorname{err}(k)\dvtx 1 \leq k\leq
M \}
.
\end{equation}
Our empirical studies reported in Section \ref{sec5} show this estimator worked
quite well.

Bickel and Levina (\citeyear{r7}, \citeyear{r8})
proposed a
method to select the bandwidth based on a repeated random splitting of
the original sample to two sub-samples of sizes $n_1$ and $n_2=n-n_1$.
Let $\hat{\Sigma}_{1}^{v}$ and $\hat{\Sigma}_{2}^{v}$ be the sample
covariances based the sub-samples of sizes $n_{1}$ and $n_{2}$
respectively, where $v$ denotes the $v$th split, for $v=1,\ldots,N$,
where $N$ is the total numbers of sample splitting.
The risk for each candidate $k$ is defined to be $R(k)={\mathrm{E}
}\Vert B_{k}(\hat
{\Sigma})-\Sigma\Vert _{(1,1)}$, where for a $p_{1}\times p_{2}$ matrix
$A=(a_{ij})$, $\Vert A\Vert _{(1,1)}=\max _{1\leq j\leq p_{2}}\sum
_{i=1}^{p_{1}}|a_{ij}|$. An empirical version of the risk is
%
\begin{equation}\label{eq:BL}
\hat{R}(k)=\frac{1}{N}\sum_{v=1}^{N}\Vert B_{k}(\hat{\Sigma
}_{1}^{v})-\hat
{\Sigma}_{2}^{v}\Vert _{(1,1)},
\end{equation}
and the bandwidth estimator is
\[
\hat{k}_{\mathrm{B L}}=\mathop{\operatorname{arg\, min}}_{0\leq k\leq p-1}\hat{R}(k).
\]
\citet{r7} recommended $n_{1}$ to be ${n}/{3}$, and the number
of random splits, $N=50$, while \citet{r8} suggested
$n_{1}=n(1-1/\log{n})$ and using the Frobenius norm instead of the
$\Vert \cdot\Vert _{(1,1)}$ norm. \citet{r27} considered a similar method
to select the bandwidth in their estimator.
We note that these approaches can be adversely impacted by high
dimensionality, due to the fact that $\hat{\Sigma}_{2}$ may be a poor
estimator of $\Sigma$ if $p$ is much larger than $n$, as found in early
works [Johnstone (\citeyear{J01}); \citet{r6}].

\section{Simulation results}\label{sec5}

In this section, we report results from simulation studies to verify
the proposed test for the bandedness and the bandwidth estimator. 
We evaluate the performance of the proposed test under several
different structures of covariance matrix for normal and gamma random vectors.
We generate $p$-dimensional independent and identical multivariate\vadjust{\goodbreak}
random vectors $X_i=(X_{i1},\ldots,X_{ip})^{\prime}$ according to a model
%
\begin{equation}\label{eq:5.1}
X_{i j}=\sum_{l=0}^{k_0}\gamma_{l}Z_{i j+l},
\end{equation}
%
where $k_0$ is the bandwidth of the covariance, $\gamma_{0}=1$ in all
settings and the other coefficients $\gamma_l$ will be specified
shortly. Two distributions are assigned to the i.i.d. $Z_{ij}$: (i) the
normal distribution $N(0,1)$; (ii) the standardized $\operatorname{Gamma}(1,0.5)$
distribution so that it has zero mean and unit variance. To mimic the
``large $p$, small~$n$'' paradigm, we choose $n=20,40,60$ and $p =
50,100,300, 600,$ respectively.

We first evaluate the size of the proposed test under the null
hypothesis $H_{k,0}\dvtx \Sigma= B_k(\Sigma)$ for $k=0$ (diagonal), $1, 2$
and $5$.
The coefficients $\gamma_l$ for $l > 0$ are: $\gamma_1 = 1$ and 0.5,
respectively, for $k= 1$; $\gamma_{1}=\gamma_{2}=1$, and $\gamma
_{1}=0.5$ and $\gamma_{2}=0.25$, respectively, for $k=2$; and $\gamma
_1= \cdots= \gamma_5 = 0.4$ for $k=5$. To assess the power, we
generate data according to (\ref{eq:5.1}) so that
$\Sigma= B_{k}(\Sigma)$ and test for $H_{k-1, 0}\dvtx \Sigma=
B_{k-1}(\Sigma)$ for $k=2$ and $5$, respectively, with the $\gamma_l$
values being the same
with those in the corresponding $k$ in the simulation for the size
reported above.
We note that this design, having the bandwidth of the null hypothesis
adjacent to the true bandwidth, is the hardest for the test, as the
null and the alternative is the closest, given the setting of the
parameters~$\{\gamma_l\}$. All the simulation results are based on 1000
simulations.

We also evaluate the test proposed in Cai and Jiang (\citeyear{r10}), based on
the asymptotic distribution of the coherence statistic $L_n$ under the
same simulation settings used for the proposed test. The test
encountered a very severe size distortion in that the real sizes are
much less than the nominal level of $5\%$, which also caused the power
of the test to be unfavorably low. For these reasons, we will not
report the simulation results of the test.
The coherence statistic is the largest Pearson correlation coefficients
among all pairs of different components in $X$, and is an extreme
value-type statistic. Extreme value statistics are known to be slowly
converging, and a computing intensive method is needed to speed up its
convergence. The asymptotic distribution established in Cai and Jiang
(\citeyear{r10}) may be the foundation to justify such a method.

\begin{table}
\caption{Empirical sizes of the proposed test at 5\% significance
for the normal and gamma random vectors generated according to model
(\protect\ref{eq:5.1})}\label{tab1}
\begin{tabular*}{\textwidth}{@{\extracolsep{\fill}}lcccccccc@{}}
\hline
& \multicolumn{4}{c}{\textbf{Normal}} & \multicolumn{4}{c@{}}{\textbf{Gamma}}\\[-5pt]
& \multicolumn{4}{c}{\hrulefill} & \multicolumn{4}{c@{}}{\hrulefill}\\
& \multicolumn{4}{c}{$\bolds{p}$} & \multicolumn{4}{c@{}}{$\bolds{p}$} \\[-5pt]
& \multicolumn{4}{c}{\hrulefill} & \multicolumn{4}{c@{}}{\hrulefill}\\
$\bolds{n}$ & \textbf{50} & \textbf{100} & \textbf{300} & \textbf{600} & \textbf{50} & \textbf{100} & \textbf{300} & \textbf{600} \\
\hline
\multicolumn{9}{@{}c@{}}{(a) $H_0\dvtx\Sigma=B_{0}(\Sigma)$}\\
20 & 0.069 & 0.065 & 0.061 & 0.066 & 0.055 & 0.056 & 0.065 & 0.075 \\
40 & 0.067 & 0.049 & 0.047 & 0.060 & 0.056 & 0.054 & 0.055 & 0.059 \\
60 & 0.066 & 0.064 & 0.045 & 0.051 & 0.068 & 0.039 & 0.065 & 0.049 \\
[6pt]
\multicolumn{9}{@{}c@{}}{$H_0:\Sigma=B_{1}(\Sigma)$}\\
\multicolumn{9}{@{}c@{}}{$\gamma_{1}=1$ } \\
20 & 0.069 & 0.061 & 0.056 & 0.060 & 0.062 & 0.058 & 0.069 & 0.069 \\
40 & 0.061 & 0.048 & 0.048 & 0.069 & 0.059 & 0.049 & 0.069 & 0.075 \\
60 & 0.045 & 0.053 & 0.056 & 0.067 & 0.048 & 0.061 & 0.068 & 0.059 \\
[3pt]
\multicolumn{9}{@{}c@{}}{$\gamma_{1}=0.5$ } \\
20 & 0.065 & 0.069 & 0.058 & 0.067 & 0.063 & 0.061 & 0.057 & 0.061 \\
40 & 0.063 & 0.052 & 0.047 & 0.068 & 0.059 & 0.055 & 0.066 & 0.071 \\
60 & 0.050 & 0.056 & 0.057 & 0.061 & 0.050 & 0.070 & 0.068 & 0.060 \\
[6pt]
\multicolumn{9}{@{}c@{}}{(c) $H_0:\Sigma=B_{2}(\Sigma)$}\\
\multicolumn{9}{@{}c@{}}{$\gamma_{1}=\gamma_{2}=1$ } \\
20 & 0.058 & 0.050 & 0.055 & 0.058 & 0.056 & 0.046 & 0.062 & 0.062 \\
40 & 0.049 & 0.042 & 0.051 & 0.058 & 0.059 & 0.048 & 0.076 & 0.071 \\
60 & 0.050 & 0.043 & 0.065 & 0.064 & 0.040 & 0.063 & 0.065 & 0.052 \\
[3pt]
\multicolumn{9}{@{}c@{}}{$\gamma_{1}=0.5$, $\gamma_{2}=0.25$ } \\
20 & 0.060 & 0.055 & 0.056 & 0.061 & 0.059 & 0.054 & 0.062 & 0.062 \\
40 & 0.055 & 0.047 & 0.055 & 0.059 & 0.058 & 0.046 & 0.071 & 0.064 \\
60 & 0.044 & 0.043 & 0.058 & 0.060 & 0.042 & 0.060 & 0.067 & 0.061 \\
[6pt]
\multicolumn{9}{@{}c@{}}{(d) $H_0:\Sigma=B_{5}(\Sigma)$ with $\gamma_{1}=\cdots =\gamma_{5}=0.4$}\\
20 & 0.045 & 0.058 & 0.067 & 0.059 & 0.050 & 0.061 & 0.054 & 0.064 \\
40 & 0.043 & 0.054 & 0.049 & 0.061 & 0.041 & 0.052 & 0.065 & 0.064 \\
60 & 0.031 & 0.046 & 0.065 & 0.069 & 0.034 & 0.040 & 0.053 & 0.048 \\
\hline
\end{tabular*}
\end{table}

Table \ref{tab1} reports the empirical sizes of the proposed test at the 5\%
nominal significance for $H_{k,0}$ with $k=0, 1, 2$ and $5$, respectively,
under both the normal and gamma distributions. Table \ref{tab2} summarizes the
empirical power of the tests whose sizes are reported in Table \ref{tab1}. To
understand the power results, Table~\ref{tab2} also contains the values of
$1-r_{k}$ for each simulation setting.
We observe from Table~\ref{tab1} that the test has reasonably empirical sizes,
around 5\%, and that the test is not sensitive to the dimensionality
indicated by its robust performance. There is some size inflation,
which is due to a~number of factors,\vadjust{\goodbreak} mainly to the dimensionality $p$,
the sample size $n$ and the approximation error of the finite sample
distribution of the test statistic by the limiting normal distribution.
We recall that the test statistic is a~linear combination of
$U$-statistics, whose convergence to the limiting normal distribution
can be slow.
In the simulations for power evaluation (reported in Table~\ref{tab2}), we
designed the simulation so that a constant $r_k$ was maintained for a
set of different $p$s, while $n$ was held fixed. The empirical powers
reported in Table~\ref{tab2} show that the power is quite reflective to the
sample size~$n$ and $1-r_{k}$, namely larger $n$ or large $1-r_{k}$
leads to higher power.
This is because as $r_{k}$ decreases, the signal of the test increases.
So it becomes easier to distinguish the null hypothesis from the alternative.
And after we controlled $n$ and $1-r_{k}$, the power was not sensitive
to $p$ at all, confirming a remark made at the end of Section \ref{sec3}.

\begin{table}
\caption{Empirical power of the proposed test at $\alpha=5\%$
for the normal and gamma random vectors generated according to model
(\protect\ref{eq:5.1})}\label{tab2}
\begin{tabular*}{\textwidth}{@{\extracolsep{\fill}}lcccccccc@{}}
\hline
& \multicolumn{4}{c}{\textbf{Normal}} & \multicolumn{4}{c@{}}{\textbf{Gamma}}\\[-5pt]
& \multicolumn{4}{c}{\hrulefill} & \multicolumn{4}{c@{}}{\hrulefill}\\
& \multicolumn{4}{c}{$\bolds{p}$} & \multicolumn{4}{c@{}}{$\bolds{p}$} \\[-5pt]
& \multicolumn{4}{c}{\hrulefill} & \multicolumn{4}{c@{}}{\hrulefill}\\
$\bolds{n}$ & \textbf{50} & \textbf{100} & \textbf{300} & \textbf{600} & \textbf{50} & \textbf{100} & \textbf{300} & \textbf{600} \\
\hline
\multicolumn{9}{@{}c@{}}{(a) $H_0\dvtx\Sigma=B_{1}(\Sigma)$ when $\Sigma=B_{2}(\Sigma)$}\\
\multicolumn{9}{@{}c@{}}{$\gamma_{1}=\gamma_{2}=1$, $1-r_{1}=1/14$ } \\
20 & 0.300 & 0.313 & 0.330 & 0.336 & 0.315 & 0.312 & 0.340 & 0.312 \\
40 & 0.683 & 0.722 & 0.711 & 0.702 & 0.710 & 0.721 & 0.752 & 0.741 \\
60 & 0.962 & 0.964 & 0.952 & 0.954 & 0.958 & 0.955 & 0.950 & 0.949 \\
[3pt]
\multicolumn{9}{@{}c@{}}{$\gamma_{1}=0.5$, $\gamma_{2}=0.25$, $1-r_{1}=1/35$ } \\
20 & 0.146 & 0.144 & 0.139 & 0.152 & 0.148 & 0.140 & 0.147 & 0.143 \\
40 & 0.269 & 0.253 & 0.258 & 0.279 & 0.256 & 0.281 & 0.311 & 0.311 \\
60 & 0.406 & 0.443 & 0.455 & 0.451 & 0.438 & 0.449 & 0.458 & 0.441 \\
[6pt]
\multicolumn{9}{@{}c@{}}{(b) $H_0\dvtx\Sigma=B_{4}(\Sigma)$ when $\Sigma=B_{5}(\Sigma)$ with}\\
\multicolumn{9}{@{}c@{}}{$\gamma_{1}=\cdots=\gamma_{5}=0.4$, $1-r_{4}=1/38.05$}\\
20 & 0.090 & 0.112 & 0.119 & 0.123 & 0.096 & 0.112 & 0.108 & 0.118 \\
40 & 0.149 & 0.181 & 0.178 & 0.200 & 0.161 & 0.169 & 0.218 & 0.196 \\
60 & 0.261 & 0.284 & 0.328 & 0.314 & 0.246 & 0.297 & 0.290 & 0.284 \\
\hline
\end{tabular*}
\end{table}

For bandwidth estimation, we generate $\{X_i\}_{i=1}^n$ according to
(\ref{eq:5.1}). While we keep $\gamma_{0}=1$, the other coefficients
$\gamma_l$ for $l > 0$ are:
\begin{longlist}[Bandwidth 10:]
\item[Bandwidth 3:] $\gamma_{i}=1$,  for  $i=1,2,3$;

\item[Bandwidth 5:] $\gamma_{i}=0.4$  for  $1\leq i\leq 5$;

\item[Bandwidth 10:] $\gamma_{i}=0.2$  for  $1\leq i\leq 5$  and  $\gamma_{i}=0.4$  for  $6\leq
i\leq10$;

\item[Bandwidth 15:] $\gamma_{i}=0.2$  for  $1\leq i\leq 10$  and  $\gamma_{i}=0.4$  for  $11\leq i\leq15$.
\end{longlist}
The covariances have bandwidth 3, 5, 10 and 15 respectively. We
evaluate two bandwidth estimators. One is $\hat{k}_{\delta, \theta}$
given in (\ref{eq:fixed}) with $\delta= 0.5$ and $\theta= 0.06$,
namely $\hat{k}_{0.5, 0.06}$, and the other is the change-point
estimator given in (\ref{eq:empirical}), applied on candidate $k$s
whose $p$-values for $H_{0 k}$ are larger than $10^{-10}$. We employ the
LOESS algorithm in R to carry our the nonparametric regression
estimation to the left of a~$k$, with a default smoothing parameter $h=0.75$.

\begin{table}
\tabcolsep=0pt
\caption{Averaged empirical bias (standard deviation)
of the five bandwidth estimators: estimator~(\protect\ref{eq:fixed}) with
$\delta=0.5$ and $\theta=0.06$ (fixed), the change-point estimator
(\protect\ref{eq:empirical}) (change-point) with $h=0.75$ and the estimators
proposed in Bickel and Levina (\protect\citeyear{r7}) (BLa), Bickel and Levina (\protect\citeyear{r8}) (BLb)
and~Rothman, Levina and Zhu (RLZ)}\label{tab3}
\begin{tabular*}{\textwidth}{@{\extracolsep{\fill}}lcld{2.9}d{2.9}d{2.9}d{2.9}@{}}
\hline
& & &\multicolumn{4}{c@{}}{\textbf{Bandwidth}}\\[-5pt]
& & &\multicolumn{4}{c@{}}{\hrulefill}\\
$\bolds{n}$ & $\bolds{p}$ & \textbf{Method} & \multicolumn{1}{c}{\textbf{3}} & \multicolumn{1}{c}{\textbf{5}} &
\multicolumn{1}{c}{\textbf{10}}& \multicolumn{1}{c@{}}{\textbf{15}} \\
\hline
20&\phantom{0}40 & Fixed & 0.58\ (1.465)& 0.07\ (0.946)&-0.5\ (1.114) &-1.63\ (1.931) \\
& &Change-point& 0.60\ (0.569)&-0.21\ (0.518)&-1.48\ (2.134)& 0.06\ (1.734) \\
&  & BLa &-0.66\ (0.855)&-0.86\ (1.287)&-4.72\ (2.202)&-9.19\ (2.246) \\
& & BLb & 0.59\ (1.036)&-0.53\ (1.460)&-3.97\ (2.932)&-6.63\ (4.403) \\
& & RLZ & 0.11\ (1.363)&-0.18\ (1.855)&-2.55\ (2.732)&-8.02\ (2.760) \\
[3pt]
& 100& Fixed & 0.14\ (0.636)& 0.1\ (0.659) &-0.22\ (0.440)&-0.96\ (0.875) \\
& &Change-point& 0.56\ (0.499)&-0.07\ (0.293)&-0.52\ (0.882)& 0.18\ (0.968) \\
 &  & BLa &-0.09\ (1.272)& 0.45\ (1.617)&-2.33\ (2.010)&-6.14\ (2.686) \\
& & BLb & 0.7 \ (1.219)&-0.26\ (1.561)&-3.88\ (2.772)&-7.29\ (3.506) \\
& & RLZ & 0.45\ (1.861)&-0.59\ (1.799)&-3.79\ (2.203)&-8.8\ (2.103) \\
[3pt]
& 200 & Fixed & 0.01\ ( 0.1 )& \multicolumn{1}{c}{0\ (0)\phantom{0000.0}} &-0.12\ (0.327)&-0.66\ (0.728) \\
& &Change-point& 0.67\ (0.473)& \multicolumn{1}{c}{0\ (0)\phantom{0000.0}} &-0.18\ (0.435)& 0.09\ (0.379) \\
& & BLa & 0.78\ (2.077)& 1.14\ (2.327)&-0.58\ (2.637)&-2.55\ (3.560) \\
& & BLb & 1.18\ (1.935)&-0.1\ (2.302) &-2.91\ (2.878)&-6.14\ (3.579) \\
& & RLZ & 0.55\ (1.641)&-0.29\ (1.719)&-4.6\ (1.928) &-9.51\ (1.823) \\
[6pt]
40& \phantom{0}80& Fixed & 0.14\ (0.551)& 0.08\ (0.464)&-0.01\ (0.1) &-0.10\ (0.302) \\
& &Change-point& 0.47\ (0.502)&-0.01\ (0.100)&-0.12\ (0.383)& 0.08\ (0.273) \\
&  & BLa &-0.24\ (0.780)& 0.23\ (1.014)&-1.32\ (1.663)&-3.55\ (2.907) \\
& & BLb & 1.5\ (1.514) & 0.94\ (1.427)& 0.06\ (2.210)&-0.17\ (3.260) \\
& & RLZ & 1.05\ (1.629)& 0.71\ (2.222)& 0.72\ (2.374)& 1.28\ (3.229) \\
[3pt]
&200 & Fixed & \multicolumn{1}{c}{0\ (0)\phantom{0000.0}} & \multicolumn{1}{c}{0\ (0)\phantom{0000.0}} & \multicolumn{1}{c}{0\ (0)\phantom{0000.0}} &-0.04\ (0.197) \\
& &Change-point& 0.55\ (0.500)& \multicolumn{1}{c}{0\ (0)\phantom{0000.0}} &-0.04\ (0.281)& 0.02\ (0.141) \\
 &  & BLa & 0.29\ (1.200)& 1.03\ (1.322)& 0.28\ (1.633)&-1.30\ (2.285) \\
& & BLb & 1.64\ (1.605)& 1.24\ (1.837)& 0.58\ (2.833)&-0.1\ (2.976) \\
& & RLZ & 1.36\ (2.435)& 1.16\ (2.465)& 2.07\ (3.647)& 1.07\ (2.861) \\
[3pt]
&400  & Fixed & \multicolumn{1}{c}{0\ (0)\phantom{0000.0}} & \multicolumn{1}{c}{0\ (0)\phantom{0000.0}} & \multicolumn{1}{c}{0\ (0)\phantom{0000.0}} & \multicolumn{1}{c@{}}{0\ (0)\phantom{0000.0}} \\
& &Change-point& 0.56\ (0.499)& \multicolumn{1}{c}{0\ (0)\phantom{0000.0}} & \multicolumn{1}{c}{0\ (0)\phantom{0000.0}} & \multicolumn{1}{c@{}}{0\ (0)\phantom{0000.0}} \\
& & BLa & 0.88\ (1.754)& 1.5\ (1.962) & 1.25\ (2.240)& 0.22\ (2.642) \\
& & BLb & 2.61\ (2.457)& 1.74\ (2.493)& 0.68\ (3.396)& 0.09\ (3.715) \\
& & RLZ & 2.19\ (2.943)& 1.98\ (3.369)& 1.17\ (3.420)&-0.39\ (2.821) \\
\hline
\end{tabular*}\vspace*{-2pt}
\end{table}

\renewcommand{\thetable}{\arabic{table}}
\setcounter{table}{2}
\begin{table}
\caption{(Continued)}
\begin{tabular*}{\textwidth}{@{\extracolsep{\fill}}lcld{1.9}d{1.9}d{1.9}d{2.9}@{}}
\hline
& & &\multicolumn{4}{c@{}}{\textbf{Bandwidth}}\\[-5pt]
& & &\multicolumn{4}{c@{}}{\hrulefill}\\
$\bolds{n}$ & $\bolds{p}$ & \textbf{Method} & \multicolumn{1}{c}{\textbf{3}} & \multicolumn{1}{c}{\textbf{5}} &
\multicolumn{1}{c}{\textbf{10}}& \multicolumn{1}{c@{}}{\textbf{15}} \\
\hline
60&120 & Fixed & 0.02\ (0.141)& 0.08\ (0.706)& 0.02\ (0.2) &-0.01\ ( 0.1 ) \\
& &Change-point& 0.52\ (0.502)& \multicolumn{1}{c}{0\ (0)\phantom{00000.0}} & \multicolumn{1}{c}{0\ (0)\phantom{00000.0}} &-0.01\ ( 0.1 ) \\
&  & BLa & 0.22\ (0.938)& 0.85\ (0.989)& 0.14\ (1.363)&-0.88\ (1.659) \\
& & BLb & 1.71\ (1.458)& 1.52\ (1.541)& 1.67\ (2.108)& 1.49\ (2.615) \\
& & RLZ & 1.24\ (1.753)& 0.71\ (1.431)& 2.03\ (2.683)& 2.13\ (2.845) \\
[3pt]
&300 & Fixed & \multicolumn{1}{c}{0\ (0)\phantom{00000.0}} & \multicolumn{1}{c}{0\ (0)\phantom{00000.0}} & \multicolumn{1}{c}{0\ (0)\phantom{00000.0}} & \multicolumn{1}{c@{}}{0\ (0)\phantom{0000.0}} \\
& &Change-point& 0.58\ (0.496)& \multicolumn{1}{c}{0\ (0)\phantom{00000.0}} & \multicolumn{1}{c}{0\ (0)\phantom{00000.0}} & \multicolumn{1}{c@{}}{0\ (0)\phantom{0000.0}} \\
 &  & BLa & 0.47\ (1.439)& 1.56\ (1.683)& 1.06\ (2.136)& 0.70\ (2.452) \\
& & BLb & 2.15\ (2.017)& 2.04\ (2.474)& 1.73\ (2.877)& 1.74\ (2.922) \\
& & RLZ & 1.68\ (2.188)& 1.02\ (2.383)& 2.45\ (3.686)& 2.75\ (3.331) \\
[3pt]
& 600 & Fixed & \multicolumn{1}{c}{0\ (0)\phantom{00000.0}} & \multicolumn{1}{c}{0\ (0)\phantom{00000.0}} & \multicolumn{1}{c}{0\ (0)\phantom{00000.0}} & \multicolumn{1}{c@{}}{0\ (0)\phantom{0000.0}} \\
& &Change-point& 0.54\ (0.501)& \multicolumn{1}{c}{0\ (0)\phantom{00000.0}} & \multicolumn{1}{c}{0\ (0)\phantom{00000.0}} & \multicolumn{1}{c@{}}{0\ (0)\phantom{0000.0}} \\
& & BLa & 1.05\ (1.702)& 1.92\ (2.102)& 2.01\ (2.393)& 1.06\ (2.490) \\
& & BLb & 3.16\ (2.631)& 2.87\ (2.699)& 2.97\ (3.532)& 1.33\ (3.254) \\
& & RLZ & 3.3\ (3.721) & 3.23\ (3.787)& 3.82\ (4.029)& 2.7\ (3.506) \\
\hline
\end{tabular*}
\end{table}

For each $\Sigma$, we compare the proposed bandwidth estimators with
the estimators advocated in Bickel and Levina (\citeyear{r7},
\citeyear{r8}) and \citet{r27}.
We choose $n$ to be 20, 40 and 60. For each $n$, $p$ is chosen 2
times, 5~times and 10 times of $n$, respectively. Following the
settings of Bickel and Levina (\citeyear{r7},
\citeyear{r8}), $n_{1}$ is chosen to be
${n}/{3}$ and $n(1-1/\log{n})$, respectively, and the number of random
splits in (\ref{eq:BL}) is $N=50$.

Table \ref{tab3} reports the average empirical bias and standard deviation of
the five bandwidth estimators based on
100 replications.
We observe from Table~\ref{tab3} that the overall performance of the proposed
estimators is better than
 those of Bickel and Levina (\citeyear{r7}, \citeyear{r8})
and \citet{r27}, with smaller standard deviation
and bias. Moreover, as $n$ is increased, both the bias and standard
deviation of the proposed estimators decreased, and are quite robust to
$p$, which is a nice property to have. For the estimators of Bickel and Levina (\citeyear{r7}, \citeyear{r8}) and \citet{r27},
the bias and the standard deviation could increase along with the
increase of $p$, and are much larger than those of the proposed estimators.
These are likely caused by the problems associated with the sample
covariance matrix when the data dimension is high.

\section{Empirical study}\label{sec6}

In this section, we report an empirical study on a~prostate cancer data
set [Adam et al. (\citeyear{r1})] from protein mass spectroscopy, which was
aimed to distinguish the healthy people from the ones with the cancer
by analyzing the constituents of the proteins in the blood.
Adam et al. (\citeyear{r1}) recorded for each blood serum sample $i$, the
intensity $X_{ij}$ for a~large number of time-of-flight values~$t_j$.
The time of flight is related to the mass over charge ratio $m/z$ of
the constituent proteins. They collected the intensity in the total of
48,538 $m/z$-sites and the full data set consisted of 157 healthy
patients and 167 with cancer.

\begin{figure}

\includegraphics{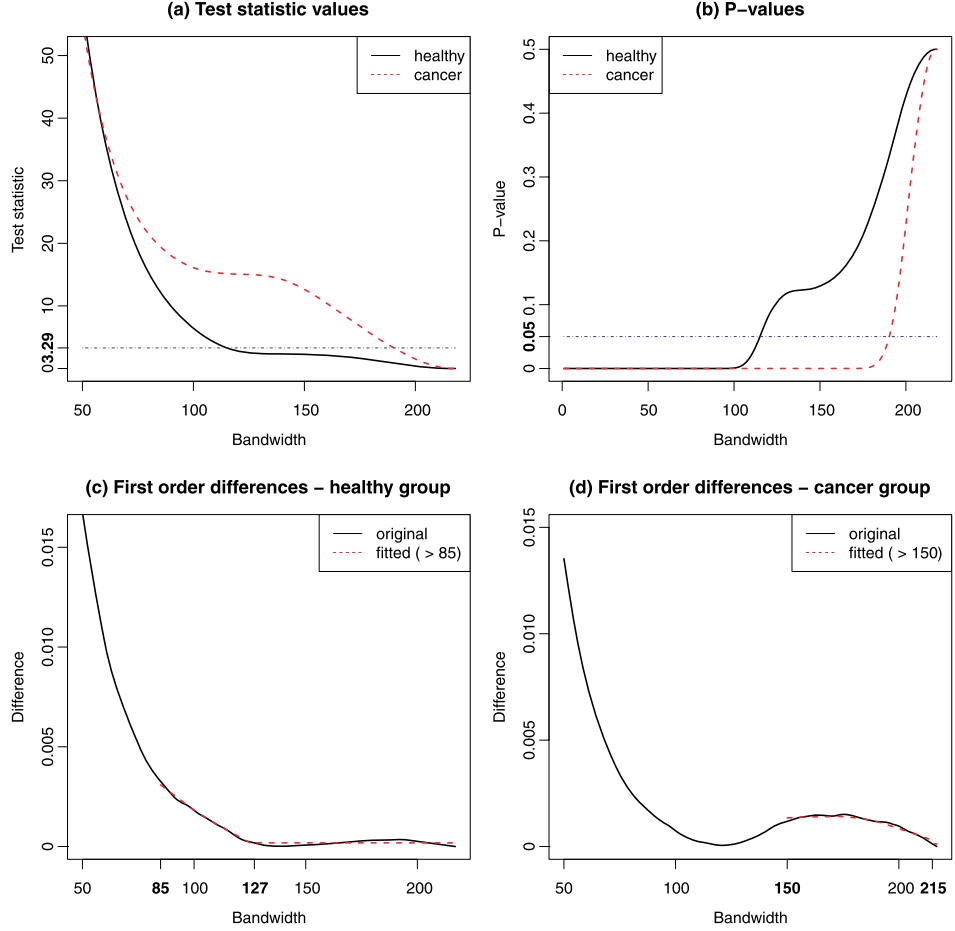}

\caption{Test statistics, $p$-values and the first order
differences $d_{nk}$ for the healthy and the cancer groups for
bandwidths larger than 50. The $p$-values of the test for $H_{0k}$ for $k
< 50$ are too small to be considered for bandwidth estimation.}\label{fig2}
\end{figure}

Tibshirani et al. (\citeyear{r29}) analyzed the data by the fused Lasso. They
ignored $m/z$-ratios below 2000 to avoid chemical artifacts, and
averaged the intensity recordings in consecutive blocks of 20. These
gave rise to a total of 2181 dimensions per observation. \citet
{r22} estimated the inverse of the covariance matrix of the intensities
by an adaptive banding approach with a nested Lasso penalty. They
carried out additional averaging of the data of Tibshirani et al.
(\citeyear{r29}) in consecutive blocks of 10, resulting in a total of 218
dimensions. 
We considered the standardized data of \citet{r22}, and tested for
the banded structure of the covariance matrix of the intensities.

The test statistics, $p$-values and the first order differences $d_{nk}$
for the healthy and cancer groups
are displayed in Figure \ref{fig2} for bandwidths $k \ge50$. We do not display
in the figure for bandwidths less than 50 since the values of the test
statistics are too large, and the associated
$p$-values for $H_{0 k}$ are too small for $k < 50$.
These bandwidth estimates together with the shapes of the curves for
the test statistics and the $p$-values in Figure~\ref{fig2} suggest that the
covariance matrix of the healthy group is likely to be banded,
while the covariance of the cancer group may not be banded at all,
given the very large bandwidth and the shape of the curve.
For the cancer group, as shown in Figure \ref{fig2}, the test statistics are
relatively flat for $120\leq k\leq140$, and then fall sharply
afterward, which indicates relatively small values in the covariance
matrix from sub-diagonal 120 to 140. However, there is a~substantial
contribution from sub-diagonals for $k > 140$.
These are echoed in the $p$-values displayed in panel (b) with almost
stationary $p$-values within the above mentioned range, followed by a
sharp increase. Panel (d) of Figure~\ref{fig2} displays
a~rather unsettled curve for $d_{nk}$, the difference between successive
statistics~$\tilde{T}_{nk}$. These are all in sharp contrasts to those
of the healthy group,
indicating rather different covariance structures between the two groups.

At $\alpha= 5\%$, we reject a $H_{k,0}$ when the statistic is larger
than 3.29. For the healthy group, the smallest $k$ such that $H_{k,0}$
is not rejected is $k=116$, while for the cancer group is 191.
We apply the bandwidth estimator (\ref{eq:fixed}) with $\delta=0.5$ and
$\theta=0.005$. The estimated bandwidth for the health group is 121 and
for the cancer group is~212.
At the same time, the bandwidth estimates, by employing Bickel and
Levina's (\citeyear{r7}) approach, are 144 for the healthy group
and 193 for the cancer group. The one for the healthy group is much
larger than the 121 we obtained earlier, using the estimator (\ref{eq:fixed}).
We then apply the proposed regression change-point bandwidth estimator
over a range of bandwidths whose associated $p$-values for testing~$H_{0
k}$ are larger than $10^{-10}$. For the healthy group, the bandwidth
range is $k\geq85$; for the cancer group the range is $k\geq150$.
We set the smoothing parameter $h=0.75$ in the LOESS procedure in R.
The regression bandwidth estimator is $\hat{k}_{h}= 127$ for the
healthy group, which is slightly larger than the 121 obtained from the
estimator (\ref{eq:fixed}). For the cancer group, the estimated
bandwidth is 215. This rather large estimated bandwidth suggests that,
compared to the healthy group, there is substantially more dependence
among the protein mass spectroscopy measurements among the cancer
patients, and, in particular, the covariance may not be banded at all
for this group of patients.

\begin{appendix}\label{app}
\section*{Appendix}

%
We first introduce some notation. For $q=0,\ldots,p$, define
\begin{eqnarray*}
B_{1,q}&=&\frac{1}{P^{2}_{n}}\sum _{l=1}^{p-q}\sum
_{i,j}^{*}(X_{il}X_{il+q})(X_{jl}X_{jl+q}),\\
B_{2,q}&=&\frac{1}{P^{3}_{n}}\sum _{l=1}^{p-q}\sum
_{i,j,k}^{*}X_{il}X_{kl+q}(X_{jl}X_{jl+q})
\end{eqnarray*}
 and
 \[
B_{3,q}=\frac{1}{P^{4}_{n}}\sum _{l=1}^{p-q}\sum
_{i,j,k,m}^{*}X_{il}X_{jl+q}X_{kl}X_{ml+q}.
\]
Then, $V_{nk}=B_{1,0}-2B_{2,0}+B_{3,0}+2\sum
_{q=1}^{k}(B_{1,q}-2B_{2,q}+B_{3,q})$, and $W_{nk}=2\sum
_{q=k+1}^{p-1}(B_{1,q}-2B_{2,q}+B_{3,q})$. Let $C_{nk}=2\sum
_{q=k+1}^{p-1}B_{1,q}$ and $U_{i}=B_{i,0}+2\sum_{q=1}^{p-1}B_{i,q}$ for
$i=1,2,3$. We first establish some lemmas for later use.

\begin{Lemma}\label{lem1}
Under Assumptions \ref{ass1} and \ref{ass2},
${\operatorname{Var}}(C_{nk})=\nu_{nk}^{2}+o\{n^{-2}{\operatorname
{tr}}^2(\Sigma^2)\}$.
\end{Lemma}

\begin{pf}
Since $C_{nk}=(P^2_n)^{-1}\sum_{i,j}^{*}\sum
_{|l_1-l_2|>k}X_{il_1}X_{il_2}X_{jl_1}X_{jl_2}$,
by the independence between different observations, we have
\[
{\mathrm{E}
}(C_{nk})=(P^2_n)^{-1}\sum_{i,j}^{*}\sum_{|l_1-l_2|>k}{\mathrm{E}
}(X_{il_1}X_{il_2}){\mathrm{E}}(X_{jl_1}X_{jl_2})
=\sum_{|l_1-l_2|>k}\sigma^2_{l_1l_2}.
\]

Note that
\[
C_{nk}^2=(P_n^2)^{-2}\sum _{i_1,j_1}^{*}\sum
_{i_2,j_2}^{*}\sum _{|l_1-l_2|>k}\sum _{|l_3-l_4|>k}
X_{i_1l_1}X_{i_1l_2}X_{i_2l_3}X_{i_2l_4}X_{j_1l_1}X_{j_1l_2}X_{j_2l_3}X_{j_2l_4}.
\]
Let $f_{l_1l_2l_3l_4}=\sum_{m}\Gamma_{l_{1}m}\Gamma_{l_{2}m}\Gamma
_{l_{3}m}\Gamma_{l_{4}m}$ and $\sigma_{l_1l_2}\sigma
_{l_3l_4}[3]=\sigma
_{l_1l_2}\sigma_{l_3l_4}+\sigma_{l_1l_3}\sigma_{l_2l_4}+\sigma
_{l_1l_4}\sigma_{l_2l_3}$. Then, ${\mathrm{E}
}(C_{nk}^{2})=(P_n^2)^{-2}(L_{n1}+L_{n2}+L_{n3})$, where
\begin{eqnarray*}
L_{n1}&=&P^4_n\sum _{|l_1-l_2|>k}\sum
_{|l_3-l_4|>k}\sigma
^2_{l_1l_2}\sigma^2_{l_3l_4},\\
L_{n2}&=&4P^3_n\sum _{|l_1-l_2|>k}\sum
_{|l_3-l_4|>k}(\Delta
f_{l_1l_2l_3l_4}+\sigma_{l_1l_2}\sigma_{l_3l_4}[3])\sigma
_{l_1l_2}\sigma
_{l_3l_4}
\end{eqnarray*}
and
\[
L_{n3}=2P^2_n\sum _{|l_1-l_2|>k}\sum
_{|l_3-l_4|>k}(\Delta
f_{l_1l_2l_3l_4}+\sigma_{l_1l_2}\sigma_{l_3l_4}[3])^2.
\]

We compute $L_{n2}$ and $L_{n3}$ part by part. First, note that
\begin{eqnarray*}
\sum _{|l_1-l_2|>k}\sum
_{|l_3-l_4|>k}f_{l_1 l_2  l_3 l_4}\sigma_{l_1l_3}\sigma_{l_2l_4}
&=&\operatorname{tr}(A^{2}\circ A^{2})-2\sum _{|l_1-l_2|\leq
k}\sum
_{l_3,l_4}f_{l_1l_2l_3l_4}\sigma_{l_1l_3}\sigma_{l_2l_4}\\
&&{}
+\sum _{|l_1-l_2|\leq k}\sum _{|l_3-l_4|\leq
k}f_{l_1l_2l_3l_4}\sigma_{l_1l_3}\sigma_{l_2l_4}.
\end{eqnarray*}
By the Cauchy--Schwarz inequality,
\[
\biggl|\sum _{|l_1-l_2|\leq k}\sum
_{l_3,l_4}f_{l_1l_2l_3l_4}\sigma_{l_1l_3}\sigma_{l_2l_4}\biggr|
\leq\operatorname{tr}^{{1}/{2}}(T^{2}){\operatorname
{tr}}^{{1}/{2}}[(\Sigma\Gamma
\circ
\Sigma\Gamma)\{(\Sigma\Gamma)'\circ(\Sigma\Gamma)'\}]
\]
  and
  \[
\biggl|\sum _{|l_1-l_2|\leq k}\sum _{|l_3-l_4|\leq
k}f_{l_1l_2l_3l_4}\sigma_{l_1l_3}\sigma_{l_2l_4}\biggr|
\leq(2k+1)^{2}\operatorname{tr}\{(\Gamma\circ\Gamma)(\Gamma
'\circ\Gamma
')(\Sigma
\circ\Sigma)\},
\]
where $T=(\Gamma\circ\Gamma)(\Gamma'\circ\Gamma')$. Note that
\[
\operatorname{tr}(T)\leq\operatorname{tr}(\Sigma^{2}),\qquad
\operatorname{tr}\{(\Gamma\circ\Gamma)(\Gamma'\circ\Gamma
')(\Sigma\circ
\Sigma)\}\leq
\operatorname{tr}(\Sigma^{4})
\]
 and
 \[
\operatorname{tr} [(\Sigma\Gamma\circ\Sigma\Gamma)\{
(\Sigma\Gamma
)'\circ(\Sigma
\Gamma)'\} ]\leq\operatorname{tr}(\Sigma^{6}).
\]
Since $\operatorname{tr}(\Sigma^{6})\leq{\operatorname
{tr}}(\Sigma^{2})\operatorname{tr}(\Sigma^{4})$,
$k=o(p^{1/4})$ and from Assumption \ref{ass2}, it follows that
\[
\sum _{|l_1-l_2|\leq k}\sum
_{l_3,l_4}f_{l_1l_2l_3l_4}\sigma_{l_1l_3}\sigma_{l_2l_4}=o\{
{\operatorname{tr}
}^2(\Sigma^2)\}
\]
 and
 \[
\sum _{|l_1-l_2|>k}\sum
_{|l_3-l_4|>k}f_{l_1l_2l_3l_4}\sigma_{l_1l_3}\sigma_{l_2l_4}=o\{
{\operatorname{tr}
}^2(\Sigma^2)\}.
\]

Similarly, it can be shown that
\begin{eqnarray*}
\sum _{|l_1-l_2|\leq k}(\Sigma^2)^{2}_{l_1l_2}&=&o\{
{\operatorname{tr}
}^2(\Sigma
^2)\},
\\
\sum _{|l_1-l_2|\leq k}\sum _{|l_3-l_4|\leq k}\sigma
^2_{l_2l_4}\sigma^2_{l_1l_3}&=&o\{\operatorname{tr}^2(\Sigma^2)\},
\\
\sum _{|l_1-l_2|\leq k}(\Sigma^2)_{l_1l_1}(\Sigma
^2)_{l_2l_2}&=&o\{
\operatorname{tr}^2(\Sigma^2)\},
\\
\sum _{|l_1-l_2|>k}\sum
_{|l_3-l_4|>k}f_{l_1 l_2 l_3 l_4}^{2}&=&o\{\operatorname{tr}^2(\Sigma
^2)\}
\end{eqnarray*}
and
\[
\sum _{  |l_1-l_2|\leq k }\sum _{
|l_3-l_4|\leq k }\sigma_{l_1l_3}\sigma_{l_2l_4}\sigma
_{l_1l_4}\sigma
_{l_2l_3}=o\{\operatorname{tr}^2(\Sigma^2)\}.
\]
 By combining these together,
\begin{eqnarray*}
{\operatorname{Var}}(C_{nk})&=&4n^{-2}\operatorname{tr}^2(\Sigma^2)+8n^{-1}
\sum _{|l_1-l_2|>k}\sum _{|l_3-l_4|>k}\sigma
_{l_1l_3}\sigma
_{l_2l_4}\sigma_{l_1l_2}\sigma_{l_3l_4}\\
&&{}+4\Delta n^{-1}\sum _{|l_1-l_2|>k}\sum
_{|l_3-l_4|>k}f_{l_1l_2l_3l_4}\sigma_{l_1l_2}\sigma_{l_3l_4}
+o (n^{-2}\operatorname{tr}^2(\Sigma^2) ).
\end{eqnarray*}
It can be checked that
\[
\sum _{|l_1-l_2|>k}\sum _{|l_3-l_4|>k}\sigma
_{l_1l_3}\sigma
_{l_2l_4}\sigma_{l_1l_2}\sigma_{l_3l_4}
=\operatorname{tr}\bigl\{\Sigma\bigl(\Sigma-B_{k}(\Sigma)\bigr)\bigr\}^{2}\vadjust{\goodbreak}
\]
and
\[
\sum _{|l_1-l_2|>k}\sum
_{|l_3-l_4|>k}f_{l_1l_2l_3l_4}\sigma_{l_1l_2}\sigma_{l_3l_4}
=\operatorname{tr}\bigl(\Gamma'\bigl(\Sigma-B_{k}(\Sigma)\bigr)\Gamma\circ
\Gamma'\bigl(\Sigma
-B_{k}(\Sigma)\bigr)\Gamma\bigr).
\]
Therefore, ${\operatorname{Var}}(C_{nk})=\nu_{nk}^{2}+o\{
n^{-2}\operatorname{tr}^2(\Sigma^2)\}$.
\end{pf}

\begin{Lemma}\label{lem2}
Under Assumptions \ref{ass1} and \ref{ass2}, for $q=0,\ldots,k$,
\[
{\operatorname{Var}}(B_{2,q})=O\{n^{-2}\operatorname{tr}^{
{1}/{2}}(\Sigma^4)\operatorname{tr}(\Sigma
^2)\}\quad
 \mbox{and}\quad
{\operatorname{Var}}(B_{3,q})=O\{n^{-4}\operatorname{tr}(\Sigma
^4)\}.
\]
\end{Lemma}

 \begin{pf}
 First consider $B_{2,q}$. Since ${\mathrm
{E}}B_{2,q}=0$ for
any $q=0,\ldots,k$, we only need to calculate ${\mathrm
{E}}B^2_{2,q}$. Note
that we can decompose $B_{2,q}^{2}$ as
\[
B_{2,q}^2=(P^{3}_{n})^{-2}\Biggl(\sum_{i=1}^{2}B_{2,q,a_{i}}+\sum
_{i=1}^{3}B_{2,q,b_{i}}+\sum_{i=1}^{2}B_{2,q,c_{i}}\Biggr),
\]
%
where
\begin{eqnarray*}
B_{2,q,a_{1}}&=&\sum _{l_1,l_2=1}^{p-q}\sum _{i,k,j1,j2}^{*}
(X_{il_1}X_{il_2})(X_{kl_1+q}X_{kl_2+q})(X_{j_1l_1}X_{j_1l_1+q})(X_{j_2l_2}X_{j_2l_2+q}),\\
B_{2,q,a_{2}}&=&\sum _{l_1,l_2=1}^{p-q}\sum _{i,k,j1,j2}^{*}
(X_{il_1}X_{il_2+q})(X_{kl_1+q}X_{kl_2})(X_{j_1l_1}X_{j_1l_1+q})(X_{j_2l_2}X_{j_2l_2+q}),\\
B_{2,q,b_{1}}&=&\sum _{l_1,l_2=1}^{p-q}\sum _{i,j,k}^{*}
(X_{il_1}X_{il_2}X_{il_2+q})(X_{jl_1}X_{jl_1+q}X_{jl_2})X_{kl_1+q}X_{kl_2+q},\\
B_{2,q,b_{2}}&=&2\sum _{l_1,l_2=1}^{p-q}\sum _{i,j,k}^{*}
(X_{il_1}X_{il_2}X_{il_2+q})(X_{jl_1}X_{jl_1+q}X_{jl_2+q})X_{kl_1+q}X_{kl_2},\\
B_{2,q,b_{3}}&=&\sum _{l_1,l_2=1}^{p-q}\sum _{i,j,k}^{*}
(X_{il_1+q}X_{il_2}X_{il_2+q})(X_{jl_1}X_{jl_1+q}X_{jl_2+q})X_{kl_1}X_{kl_2},\\
B_{2,q,c_{1}}&=&\sum _{l_1,l_2=1}^{p-q}\sum _{i,j,k}^{*}
(X_{il_1}X_{il_2})(X_{kl_1+q}X_{kl_2+q})(X_{jl_1}X_{jl_1+q}X_{jl_2}X_{jl_2+q})
\end{eqnarray*}
 and
 \[
B_{2,q,c_{2}}=\sum _{l_1,l_2=1}^{p-q}\sum _{i,j,k}^{*}
(X_{il_1}X_{il_2+q})(X_{kl_1+q}X_{kl_2})(X_{jl_1}X_{jl_1+q}X_{jl_2}X_{jl_2+q}).
\]
We need to show that the expectations of all the terms above are controlled by the order $n^4\operatorname{tr}^{1/2}(\Sigma^4)\operatorname{tr}(\Sigma^2)$.
First, note that
${\mathrm{E}}(B_{2,q,a_{1}})=\break P^{4}_{n}\sum _{l_1,l_2=1}^{p-q}
\sigma_{l_1l_2}\sigma_{l_1+q
l_2+q}\*\sigma_{l_1l_1+q}\sigma_{l_2l_2+q}$.\vspace*{2pt}\vadjust{\goodbreak}
By the Cauchy--Schwarz inequality, it can be shown that
\[
|{\mathrm{E}}(B_{2,q,a_{1}})|
=P^{4}_{n}O(\operatorname{tr}^{{1}/{2}}(\Sigma
^{4})\operatorname{tr}(\Sigma^{2})).
\]

Employing a similar derivation, we can show that the same result holds
for all the other terms, which lead to the first part of Lemma \ref{lem2}. The
second part can be proved following the same track.
\end{pf}

\begin{Lemma}\label{lem3}
Under Assumptions \ref{ass1} and \ref{ass2},
${\operatorname{Var}}(U_{i})=o\{n^{-2}\operatorname{tr}^2(\Sigma
^2)\}$ for \mbox{$i=2,3$}.
\end{Lemma}

\begin{pf}
The proof is similar to Lemma \ref{lem2}.
\end{pf}

\begin{Lemma}\label{lem4}
$\!\!\!$Under Assumptions \ref{ass1} and \ref{ass2}, ${\operatorname
{Var}} (\sum
_{q=k+1}^{p-1}B_{i,q} )\,{=}\,o\{n^{-2}\operatorname{tr}^2(\Sigma^2)\}
$ for $i=2,3$.
\end{Lemma}

\begin{pf}
Noting that $\sum_{q=k+1}^{p}B_{i,q}=U_{i}-\sum
_{q=1}^{k}B_{i,q}$, the lemma follows by applying Lemmas \ref{lem2}, \ref{lem3},
$k=o(p^{{1}/{4}})$ and Assumption \ref{ass2}.
\end{pf}

In the following, we provide the proof of Propositions \ref{pro1} and \ref{pro2}.

\begin{pf*}{Proof of Proposition \ref{pro1}}
Rewrite $W_{nk}$ as
\[
W_{nk}=C_{nk}-2\sum_{q=k+1}^{p}B_{2,q}+\sum_{q=k+1}^{p}B_{3,q}.
\]
Since ${\mathrm{E}}(C_{nk})=\sum_{|i-j|>k}\sigma
_{ij}^{2}=\operatorname{tr}[\{\Sigma-
B_k(\Sigma)\}^2]$ and ${\mathrm{E}}(B_{i,q})=0$ for $i=2,3$ and any
$q=0,1,\ldots,p-1$, the first statement is readily obtained. The second
statement follows by applying Lemmas \ref{lem1}, \ref{lem4} and the fact that $\nu
_{n k}^2 \ge4n^{-2}\operatorname{tr}^2(\Sigma^2)$.
\end{pf*}

\begin{pf*}{Proof of Proposition \ref{pro2}}
It can be carried out following the same routes as those in Lemmas \ref{lem1} and
\ref{lem2}. Specifically, it can be shown that ${\operatorname{Var}}(V_{nk})=O\{
a_{np}{\operatorname{tr}
}^{2}(\Sigma^{2})\}$. Hence,
${\operatorname{Var}}\{V_{nk}/\operatorname{tr}(\Sigma^{2})\}
=O(a_{np})\to0$.
\end{pf*}

It is clear from the proof of Proposition \ref{pro1} that
$W_{nk}=C_{nk}+o_{p}(\nu_{nk})$. Therefore, in order to derive the
asymptotical distribution of the statistic, we only need to consider
the asymptotical normality of $C_{nk}$.
Let $\mathscr{F}_0=\{\varnothing,\Omega\}$, and $\mathscr
{F}_t=\sigma\{
X_1,\ldots,X_t\}$ for $t=1,2,\ldots,n$, be a sequence of $\sigma$-field
generated by the data sequence. Let ${\mathrm{E}}_t(\cdot)$ denote the
conditional expectation with respect to $\mathscr{F}_t$. Write
$C_{nk}-{\mathrm{E}}(C_{nk})=\sum^n_{t=1}D_{tk}$, where
$D_{tk}=({\mathrm{E}}_t-{\mathrm{E}
}_{t-1})C_{nk}$. Then for every~$n$, $D_{tk},1\leq t\leq n$, is a
martingale difference sequence with respect to the $\sigma$-fields $\{
\mathscr{F}_t\}_{t=0}^{\infty}$.

\begin{Lemma}\label{lem5}
Let $\sigma^2_{tk}={\mathrm{E}}_{t-1}(D^2_{tk})$. Under
Assumptions \ref{ass1} and \ref{ass2}, as $n\rightarrow\infty$,
%
\begin{equation}\label{eA1}
\frac{\sum^n_{t=1}\sigma^2_{tk}}{{\operatorname
{Var}}(C_{nk})}\stackrel{p}{\rightarrow}
1\quad
\mbox{and}\quad \frac{\sum^n_{t=1}{\mathrm
{E}}(D^4_{tk})}{{\operatorname{Var}
}^2(C_{nk})}\rightarrow0.\vadjust{\goodbreak}
\end{equation}
\end{Lemma}

 \begin{pf}
$\!\!\!\!$We first establish the first part of (\ref{eA1}). Noting
that ${\mathrm{E}}(\sum^n_{t=1}\sigma^2_{tk})={\operatorname
{Var}}(C_{nk})$, we need only to
show ${\operatorname{Var}}(\sum^n_{t=1}\sigma
^2_{tk})=o({\operatorname{Var}}^2(C_{nk}))$. Note that
\begin{eqnarray*}
D_{tk}&=&\frac{2}{n(n-1)} \Biggl[\sum
_{|l_1-l_2|>k}(X_{tl_1}X_{tl_2}-\sigma
_{l_1l_2}) \Biggl\{\sum _{i=1}^{t-1}(X_{il_1}X_{il_2}
-\sigma_{l_1l_2}) \Biggr\} \Biggr]\\
&&{}+\frac{2}{n} \biggl(\sum_{|l_1-l_2|>k}X_{tl_1}X_{tl_2}\sigma
_{l_1l_2}-\sum
_{|l_1-l_2|>k}\sigma^2_{l_1l_2} \biggr).
\end{eqnarray*}
Denote $Q^{l_1l_2}_{t-1}=\sum^{t-1}_{i=1}(X_{il_1}X_{il_2}-\sigma
_{l_1l_2})$. Let $Q_{t-1}$ be the matrix with the $(l_1,l_2)$th entry being
$Q^{l_1l_2}_{t-1}$ and $M_{t-1}=\Gamma'Q_{t-1}\Gamma$; then
%
\[
\sum^n_{t=1}\sigma^2_{tk}=\sum _{i=1}^{3}R_{1i}+\Delta\sum
_{i=1}^{3}R_{2i}
+\sum _{i=1}^{4}R_{3i}+\Delta\sum
_{i=1}^{4}R_{4i}+n\gamma,
\]
where $\gamma$ is a constant and
%
\begin{eqnarray*}
R_{11}&=&\frac{4}{n^2(n-1)^2}\sum _{t=1}^{n}{\operatorname
{tr}}(M_{t-1}^2),\\
R_{12}&=&-\frac{8}{n^2(n-1)^2}\sum _{t=1}^{n}\sum
_{|l_1-l_2|\leq
k}Q_{t-1}^{l_1l_2}(\Sigma Q_{t-1}\Sigma)_{l_1l_2},\\
R_{13}&=&\frac{4}{n^2(n-1)^2}\sum _{t=1}^{n}\sum
_{|l_1-l_2|\leq
k}\sum_{|l_3-l_4|\leq k}Q_{t-1}^{l_1l_2}Q_{t-1}^{l_3l_4}\sigma
_{l_1l_3}\sigma_{l_2l_4},\\
R_{21}&=&\frac{4}{n^2(n-1)^2}\sum _{t=1}^{n}{\operatorname
{tr}}(M_{t-1}\circ
M_{t-1}),\\
R_{22}&=&-\frac{8}{n^2(n-1)^2}\sum _{t=1}^{n}\sum
_{m}\sum
_{|l_1-l_2|\leq k}Q^{l_1l_2}_{t-1}M^{mm}_{t-1}\Gamma_{l_1m}\Gamma
_{l_2m},\\
R_{23}&=&\frac{4}{n^2(n-1)^2}\sum _{t=1}^{n}\sum
_{m}\sum
_{|l_1-l_2|\leq k}\sum_{|l_3-l_4|\leq k}Q_{t-1}^{l_1l_2}Q_{t-1}^{l_3l_4}
\Gamma_{l_1m}\Gamma_{l_2m}\Gamma_{l_3m}\Gamma_{l_4m},\\
R_{31}&=&\frac{8}{n^2(n-1)}\sum _{t=1}^{n}{\operatorname
{tr}}(\Sigma
Q_{t-1}\Sigma^2),\\
R_{32}&=&-\frac{8}{n^2(n-1)}\sum _{t=1}^{n}\sum_{|l_1-l_2|\leq
k}Q_{t-1}^{l_1l_2}(\Sigma^3)_{l_1l_2},
\\
R_{33}&=&-\frac{8}{n^2(n-1)}\sum _{t=1}^{n}\sum_{|l_1-l_2|\leq
k}(\Sigma Q_{t-1}\Sigma)_{l_1l_2}\sigma_{l_1l_2},\\
R_{34}&=&\frac{8}{n^2(n-1)}\sum _{t=1}^{n}\sum_{|l_1-l_2|\leq
k}\sum_{|l_3-l_4|\leq k}Q_{t-1}^{l_1l_2}\sigma_{l_3l_4}\sigma
_{l_1l_3}\sigma_{l_2l_4},\\
R_{41}&=&\frac{8}{n^2(n-1)}\sum _{t=1}^{n}{\operatorname
{tr}}(M_{t-1}\circ
A^2),\\
R_{42}&=&-\frac{8}{n^2(n-1)}\sum _{t=1}^{n}\sum
_{m}\sum
_{|l_1-l_2|\leq k}Q^{l_1l_2}_{t-1}\Gamma_{l_1m}\Gamma
_{l_2m}(A^2)_{mm},\\
R_{43}&=&-\frac{8}{n^2(n-1)}\sum _{t=1}^{n}\sum
_{m}\sum
_{|l_1-l_2|\leq k}\sigma_{l_1l_2}\Gamma_{l_1m}\Gamma
_{l_2m}M^{mm}_{t-1}
\end{eqnarray*}
and
\[
R_{44}=\frac{8}{n^2(n-1)}\sum _{t=1}^{n}\sum _{m}\sum
_{|l_1-l_2|\leq k}\sum_{|l_3-l_4|\leq k}
Q_{t-1}^{l_1l_2}\sigma_{l_3l_4}\Gamma_{l_1m}\Gamma_{l_2m}\Gamma
_{l_3m}\Gamma_{l_4m}.
\]

To prove ${\operatorname{Var}}(\sum^n_{t=1}\sigma
^2_{tk})=o({\operatorname{Var}}^2(C_{nk}))$, we intend
to prove the variance of each $R_{ij}$ is of small order of
$n^{-4}{\operatorname{tr}
}^{4}(\Sigma^{2})$.

For $R_{12}$, denote for any $1\leq i,j\leq n$,
\[
Y^{12}_{ij}=\sum_{|l_1-l_2|\leq k}(X_{il_1}X_{il_2}-\sigma_{l_1l_2})\{
(\Sigma X_{j}X'_{j}\Sigma)_{l_1l_2}-(\Sigma^3)_{l_1l_2}\}.
\]
Then $\sum_{|l_1-l_2|\leq k}Q_{t-1}^{l_1l_2}(\Sigma Q_{t-1}\Sigma
)_{l_1l_2}=\sum^{t-1}_{i=1}Y^{12}_{ii}+\sum^{t-1}_{i\neq j}Y^{12}_{ij}$.
Note that\break
${\mathrm{E}}Y^{12}_{ij}=0$ for any $i\neq j$ and ${\mathrm{E}
}(Y^{12}_{i_1j_1}Y^{12}_{i_2j_2})=0$ for any $(i_1,i_2,j_1,j_2)$,
except $\{i_1=i_2,j_1=j_2\}$ and $\{i_1=j_1,i_2=j_2\}$. Thus for any $t<l$,
\begin{eqnarray*}
&&{\operatorname{Cov}}\biggl(\sum_{|l_1-l_2|\leq k}Q_{t-1}^{l_1l_2}(\Sigma
Q_{t-1}\Sigma
)_{l_1l_2},\sum_{|l_1-l_2|\leq k}Q_{l-1}^{l_1l_2}(\Sigma Q_{l-1}\Sigma
)_{l_1l_2}\biggr)\\
&&\qquad=(t-1){\operatorname{Var}}(Y^{12}_{11})+(t-1)(t-2){\operatorname
{Var}}(Y^{12}_{12}).
\end{eqnarray*}
%
We only need to verify that ${\operatorname{Var}}(y_{11}^{12})$ and
${\operatorname{Var}}(y_{12}^{12})$
are of small orders of $\operatorname{tr}^{4}(\Sigma^{2})$. Note that
%
\begin{eqnarray*}
{\mathrm{E}}(Y^{12}_{11})^2&=&{\mathrm{E}}\sum_{|l_1-l_2|\leq k}\sum
_{|l_3-l_4|\leq
k}(X_{1l_1}X_{1l_2}-\sigma_{l_1l_2})(X_{1l_3}X_{1l_4}-\sigma
_{l_3l_4})\\
&&\hspace*{74pt}{}\times\{(\Sigma X_{1}X'_{1}\Sigma)_{l_1l_2}-(\Sigma^3)_{l_1l_2}\}\\
&&\hspace*{74pt}{}\times\{
(\Sigma X_{1}X'_{1}\Sigma)_{l_3l_4}-(\Sigma^3)_{l_3l_4}\}\\
&\leq& \gamma_{12}\sum_{|l_1-l_2|\leq k}\sum_{|l_3-l_4|\leq
k}(\sigma
_{l_1l_2}^{2}+\sigma_{l_1l_1}\sigma_{l_2l_2})^{{1}/{2}}
(\sigma_{l_3l_4}^{2}+\sigma_{l_3l_3}\sigma_{l_4l_4})^{{1}/{2}}\\
&&\hspace*{81pt}{}\times\{(\Sigma^{3})_{l_1l_2}^{2}+(\Sigma^{3})_{l_1l_1}(\Sigma
^{3})_{l_2l_2}\}^{{1}/{2}}\\
&&\hspace*{81pt}{}\times
\{(\Sigma^{3})_{l_3l_4}^{2}+(\Sigma^{3})_{l_3l_3}(\Sigma
^{3})_{l_4l_4}\}
^{{1}/{2}}\\
&\leq&\gamma_{12}\sum_{|l_1-l_2|\leq k}(\sigma_{l_1l_2}^{2}+\sigma
_{l_1l_1}\sigma_{l_2l_2})
\sum_{|l_1-l_2|\leq k}\{(\Sigma^{3})_{l_1l_2}^{2}+(\Sigma
^{3})_{l_1l_1}(\Sigma^{3})_{l_2l_2}\}\\
&\leq&\gamma_{12}(2k+1)^{2}\operatorname{tr}(\Sigma
^{2})\operatorname{tr}(\Sigma^{6}),
\end{eqnarray*}
where $\gamma_{12}$ is a constant. Since $\operatorname{tr}(\Sigma
^{6})\leq{\operatorname{tr}
}^{{3}/{2}}(\Sigma^{4})$,
\begin{eqnarray*}
(2k+1)^{2}\operatorname{tr}(\Sigma^{2})\operatorname{tr}(\Sigma
^{6})&=&O\{k^{2}\operatorname{tr}(\Sigma
^{2})\operatorname{tr}^{{3}/{2}}(\Sigma^{4})\}\\
&=&O\{k^{2}p^{-{3}/{2}}\operatorname{tr}^{4}(\Sigma^{2})\}\\
&=&o\{
{\operatorname{tr}
}^{4}(\Sigma
^{2})\},
\end{eqnarray*}
which indicates that ${\operatorname{Var}}(Y_{11}^{12})=o\{
\operatorname{tr}^{4}(\Sigma^{2})\}$.
Similarly, we can also show that ${\operatorname{Var}}(Y^{12}_{12})=o\{
\operatorname{tr}^4(\Sigma
^2)\}$.
Thus
\begin{eqnarray*}
{\operatorname{Var}}(R_{12})&=&\frac{64}{n^4(n-1)^4}{\operatorname
{Var}} \Biggl\{\sum ^{n}_{t=1}\sum
_{|l_1-l_2|\leq k}Q_{t-1}^{l_1l_2}(\Sigma Q_{t-1}\Sigma)_{l_1l_2}
\Biggr\}\\
&=&o\{n^{-4}\operatorname{tr}^4(\Sigma^2)\}.
\end{eqnarray*}

Following the same procedure, we can prove that for all the other
$R_{ij}$, ${\operatorname{Var}}(R_{ij})=o\{n^{-4}{\operatorname
{tr}}^{4}(\Sigma^2)\}$. Since ${\operatorname{Var}
}^2(C_{nk})\geq n^{-4}\operatorname{tr}^4(\Sigma^2)$, we have\break
${\operatorname{Var}}(R_{ij})=o\{
{\operatorname{Var}
}^2(C_{nk})\}$. Thus we have ${\operatorname{Var}}(\sum^n_{t=1}\sigma
^2_{tk})=o({\operatorname{Var}
}^2(C_{nk}))$, and hence the first part of (\ref{eA1}).

For the second part of (\ref{eA1}), by simple algebra, we can rewrite
$D_{tk}$ as $D_{tk}=S_{t1}-S_{t2}+S_{t3}-S_{t4}$, where
%
\begin{eqnarray*}
S_{t1}&=&\frac{2}{n(n-1)} \{X'_{t}Q_{t-1}X_{t}-{\operatorname
{tr}}(Q_{t-1}\Sigma
) \},\\
S_{t2}&=&\frac{2}{n(n-1)} [X'_{t}B_{k}(Q_{t-1})X_{t}-{\operatorname
{tr}}\{
B_{k}(Q_{t-1})\Sigma\} ],\\
S_{t3}&=&\frac{2}{n} \{X'_{t}\Sigma X_{t}-\operatorname{tr}(\Sigma
^2) \}
\end{eqnarray*}
and
\[
S_{t4}=\frac{2}{n} [X'_{t}B_{k}(\Sigma)X_{t}-\operatorname{tr}\{
B_{k}(\Sigma
)\Sigma\} ].
\]
Since $D^4_{tk}\leq\tilde{\gamma
}(S_{t1}^4+S_{t2}^4+S_{t3}^4+S_{t4}^4)$, we have for a positive
constant $\tilde{\gamma}$,
\[
\sum _{t=1}^{n}{\mathrm{E}}(D^4_{tk})\leq\tilde{\gamma}
\Biggl\{\sum _{t=1}^{n}{\mathrm{E}}(S_{t1}^4)+\sum
_{t=1}^{n}{\mathrm{E}
}(S_{t2}^4)+\sum _{t=1}^{n}{\mathrm{E}}(S_{t3}^4)
+\sum _{t=1}^{n}{\mathrm{E}}(S_{t4}^4) \Biggr\}.
\]
In the following, we will prove the four terms on the right are of
small orders of ${\operatorname{Var}}^{2}(C_{nk})$, respectively. To
this end, note that
\[
{\mathrm{E}} \{X'_{t}Q_{t-1}X_{t}-\operatorname{tr}(Q_{t-1}\Sigma)
\}^4
\leq\tilde{\gamma_{1}}{\mathrm{E}}\{{\operatorname
{tr}}^2(M^2_{t-1})\},
\]
%
where $\tilde{\gamma_{1}}$ is a positive constant. Since ${\mathrm
{E}}\{{\operatorname{tr}
}(M^2_{t-1})\}=(t-1)O\{\operatorname{tr}^2(\Sigma^2)\}$, and
${\operatorname{Var}}\{{\operatorname{tr}
}(M^2_{t-1})\}=t^2O(\operatorname{tr}^2(\Sigma^2){\operatorname
{tr}}(\Sigma^4))$, then we have
${\mathrm{E}}\{\operatorname{tr}^2(M^2_{t-1})\}=t^2\times O\{
\operatorname{tr}^4(\Sigma^2)\}$. Thus,
%
\begin{eqnarray*}
\sum _{t=1}^{n}{\mathrm{E}}(S_{t1}^4)&=&\frac
{16}{n^4(n-1)^4}\sum
 _{t=1}^{n}
{\mathrm{E}} \{X'_{t}Q_{t-1}X_{t}-\operatorname{tr}(Q_{t-1}\Sigma)
\}^4\\
&\leq&\frac{16}{n^4(n-1)^4}\sum _{t=1}^{n}t^2O\{{\operatorname{tr}
}^4(\Sigma^2)\}
=\frac{1}{n^5}O\{\operatorname{tr}^4(\Sigma^2)\}=o\{{\operatorname
{Var}}^2(C_n)\}.
\end{eqnarray*}

Similarly, we can show that for $i=2,3$ and $4$, $\sum
_{t=1}^{n}{\mathrm{E}
}(S_{ti}^4)=o\{{\operatorname{Var}}^2(C_n)\}$. Combining all the four
parts together,
we have $\sum_{t=1}^{n}{\mathrm{E}}(D^4_{k,t})=o\{{\operatorname
{Var}}^2(C)\}$, which leads to
the second part of (\ref{eA1}).
\end{pf}

Denote $I_{nk}=\{W_{nk}-{\mathrm{E}}(W_{nk})\}/V_{nk}$ and
$J_{nk}={\mathrm{E}
}(W_{nk})/V_{nk}$. Then $\tilde{T}_{nk}=I_{nk}+J_{nk}$. For $k_{0}$
diverging, but satisfying (\ref{eq:k0}), we intend to prove $n^{\delta
}(J_{nk}-J_{n k+1})$ diverging to $\infty$ uniformly on $k<k_{0}$ for
any $\delta>0$. And $n^{\delta}I_{nk}$ uniformly converges to $0$ in
probability for any $\delta\leq1/2$ and $k\leq M$, where $M>k_{0}$ and
$M=o(p^{1/4})$.

\begin{Lemma}\label{lem6}
Under Assumptions \ref{ass1}, \ref{ass2} and (\ref{eq:k0}),
if $\liminf_n \{\inf_{k<k_{0}}(r_{k+1}-r_{k})\}>0$ and $\{\sigma
_{ll}\}
_{l=1}^{p}$ are uniformly bounded away from $0$ and $\infty$, for any
$\delta\leq0.5$, as $n\to\infty$:
\begin{longlist}[(b)]
\item[(a)] $P (n^{\delta}(J_{nk}-J_{n k+1})>\xi$, for any $k<
k_{0} )\to1$ for any $\xi>0$;

\item[(b)] $P (n^{\delta}|I_{nk}| \leq\varepsilon$,   for any
$k\leq k_{0} )\to1$ for any $0<\varepsilon<1$;

\item[(c)] $P (n^{\delta}|I_{nk}| \leq\varepsilon$,   for any
$k_{0}<k\leq M )\to1$ for any $0<\varepsilon<1$, where $k_{0}<M$
and $M=o(p^{1/4})$.
\end{longlist}
\end{Lemma}

\begin{pf}  (a) If $\{\sigma_{ll}\}_{l=1}^{p}$ is bounded away from
$\infty$, similarly to the proof of Lemmas \ref{lem1} and \ref{lem2}, it can be
checked that ${\operatorname{Var}}(V_{n k})=O(k^{2}{\operatorname
{tr}}(\Sigma^{2})/n)$. Therefore, by
Chebyshev's inequality, for any $\varepsilon>0$,
\[
P \biggl( \biggl|\frac{V_{n k}-{\mathrm{E}}(V_{n k})}{{\operatorname
{tr}}(\Sigma^{2})}
\biggr|>\varepsilon r_{k}^{2} \biggr)
\leq\frac{{\operatorname{Var}}(V_{n k})}{\varepsilon
^{2}\operatorname{tr}^{2}(\Sigma
^{2})r_{k}^{4}}\leq\frac{Ck^{2}}{\varepsilon^{2}npr_{k}^{3}}
\leq\frac{Ck^{2}k_{0}^{3}}{\varepsilon^{2}np},
\]
where the last inequality comes from the fact that $r_{k}^{-1}\leq
2k_{0}+1$. Hence,
\[
P \biggl(\max_{0 \leq k\leq k_{0}}\biggl |\frac{V_{n k}-{\mathrm{E}}(V_{n
k})}{{\operatorname{tr}
}(\Sigma^{2})r_{k}^{2}} \biggr|\leq\varepsilon \biggr)
\geq1- \sum_{k=0}^{k_{0}}\frac{Ck^{2}k_{0}^{3}}{\varepsilon^{2}np}
\geq1- \frac{C k_{0}^{6}}{\varepsilon^{2}np},
\]
which converge to $1$ since $k_{0}$ satisfies (\ref{eq:k0}).
Consider $\varepsilon<1/2$, and denote
\[
\Omega= \{\omega\dvtx|V_{n k}-{\mathrm{E}}(V_{n k})|\leq\varepsilon
r_{k}^{2}
\operatorname{tr}(\Sigma^{2}), \mbox{  for any $k\leq k_{0}$}
\}.\vadjust{\goodbreak}
\]
By the above argument, $P(\Omega)\to1$ as $n\to\infty$. For any
$\omega\in\Omega$, we have
\[
1-\varepsilon r_{k}\leq1/(1+\varepsilon r_{k})\leq{\operatorname
{tr}}[\{
B_{k}(\Sigma
)\}^2]/V_{nk}\leq1/(1-\varepsilon r_{k})\leq1+2\varepsilon r_{k}
\]
for any $k< k_{0}$. Hence, for any $\omega\in\Omega$,
\begin{eqnarray*}
n^{\delta}(J_{nk}-J_{n k+1})&\geq& n^{\delta}(r_{k+1}-r_{k})+n^{\delta}(
\varepsilon r_{k}+2\varepsilon r_{k+1}-3\varepsilon)\\
&\geq& n^{\delta}(r_{k+1}-r_{k})-3n^{\delta}\varepsilon,
\end{eqnarray*}
which implies that $n^{\delta}(J_{nk}-J_{n k+1})$ diverge uniformly on
$k<k_{0}$, by choosing~$\varepsilon$ small enough. Therefore, for any
$\xi>0$, by choosing $\varepsilon$ small enough, there exists a $N>0$
such that for any $n>N$,
\[
P \bigl(n^{\delta}(J_{nk}-J_{n k+1})>\xi\mbox{  for any $k<
k_{0}$}
\bigr)\geq P(\Omega).
\]
The conclusion follows by noting that $P(\Omega)\to1$ as $n\to\infty$.
The other two parts of the conclusion can be obtained similarly. For
simplicity in the presentation, we omit them here.
\end{pf}

\begin{pf*}{Proof of Theorem \ref{teo1}}
By Lemmas \ref{lem1}, \ref{lem5} and the martingale
central limit theorem [\citet{r9}], it is readily shown that as
$n\to\infty$,
\[
\frac{C_{nk}-{\mathrm{E}}(C_{nk})}{\nu_{nk}}\stackrel{D}{\rightarrow} N(0,1).
\]
Substituting $C_{nk}$ for $W_{nk}$, Theorem \ref{teo1} follows by noting $W_{nk}=C_{nk}+o_p(\nu_{nk})$.
\end{pf*}

\begin{pf*}{Proof of Theorem \ref{teo2}}
Note that ${\operatorname{Var}}\{V_{nk}/\operatorname{tr}(\Sigma
^{2})\}\rightarrow0$, ${\mathrm{E}
}\{
V_{nk}/\break \operatorname{tr}(\Sigma^{2})\}= r_{k}$ and $\limsup
_n r_{k}\leq1$.
It can be shown that for any $\eta>0$,\break $\lim_{n\to\infty
}P(B_{n,\eta})=1$
where $B_{n,\eta}=\{V_{nk}<(1+\eta)\operatorname{tr}(\Sigma^{2})\}
$. This means
that for any $\varepsilon>0$, there exists a positive integer $N$, such
that for all $n>N$, $P(B_{n,\eta})>1-\varepsilon$. Then from (\ref{eq:beta}),
\begin{eqnarray*}
\beta_{nk}&\geq& P \biggl(\frac{W_{nk}-\operatorname{tr}(\Sigma
^2)+\operatorname{tr}[\{
B_{k}(\Sigma
)\}^{2}]}{\nu_{nk}}\geq z_{\alpha}\frac{V_{nk}}{{\operatorname
{tr}}(\Sigma
^{2})}-\delta_{nk}, B_{n,\eta} \biggr)\\[-2pt]
&\geq& P \biggl(\frac{W_{nk}-\operatorname{tr}(\Sigma^2)+{\operatorname
{tr}}[\{B_{k}(\Sigma)\}
^{2}]}{\nu_{nk}}\geq z_{\alpha}(1+\eta)-\delta_{nk}, B_{n,\eta} \biggr)\\[-2pt]
&\geq& P \biggl(\frac{W_{nk}-\operatorname{tr}(\Sigma^2)+{\operatorname
{tr}}[\{B_{k}(\Sigma)\}
^{2}]}{\nu_{nk}}\geq z_{\alpha}(1+\eta)-\delta_{nk} \biggr)-P(B_{n,\eta}^{c}).
\end{eqnarray*}
Therefore, from Theorem \ref{teo1},
\begin{eqnarray*}
\liminf_{n\to\infty}\beta_{nk}&\geq& \liminf_{n\to\infty}P \biggl\{
\frac
{W_{nk}-\operatorname{tr}(\Sigma^2)+\operatorname{tr}[\{
B_{k}(\Sigma)\}^{2}]}{\nu
_{nk}}\geq
z_{\alpha}(1+\eta)-\delta_{nk} \biggr\}\\[-2pt]
&&{}-\limsup_{n\to\infty}P(B_{n,\eta}^{c})\\[-2pt]
&\geq& 1-\Phi\Bigl\{z_{\alpha}(1+\eta)-\liminf_{n\to\infty}\delta
_{nk}\Bigr\}
-\varepsilon.
\end{eqnarray*}
The first part of the theorem follows by taking $\varepsilon\to0$ and
$\eta\to0$.\vadjust{\goodbreak}

(ii) The condition $a_{np}^{-1/2}(1-r_{k})\to\infty$ implies that
$\delta_{nk}\to\infty$ as $n\to\infty$. Hence, $\beta_{nk}\to1$.
\end{pf*}

\begin{pf*}{Proof of Theorem \ref{teo3}}
First consider the case where $k_{0}$ is bounded. Consider $M$ to be a
fixed sufficiently large integer. Recall that $\tilde
{T}_{nk}=I_{nk}+J_{nk}$, where
\[
I_{nk}=\{W_{nk}-{\mathrm{E}}(W_{nk})\}/V_{nk}\quad \mbox{and}\quad
J_{nk}={\mathrm{E}
}(W_{nk})/V_{nk}.
\]
By (\ref{eq:Wnk}), since $a_{np}^{{1}/{2}}=O(n^{-1})$, we have
$n^{\delta}I_{nk}=O_{p}(n^{\delta}a_{np}^{{1}/{2}})\to0$, for any
$k\leq M$. Note that
\[
n^{\delta}(r_{k}^{-1}-r_{k+1}^{-1})=n^{\delta}\frac
{r_{k+1}-r_{k}}{r_{k+1}r_{k}}\geq n^{\delta}(r_{k+1}-r_{k}).
\]
Thus, from (\ref{eq:Wnk}), for $k<k_{0}$, the condition $\liminf
_n (r_{k+1}-r_{k})>0$ implies that $n^{\delta}(J_{nk}-J_{n k+1})\sim
n^{\delta}\to\infty$ in probability, where $\delta\in(0,1)$. Therefore,
${d}_{n k}^{(\delta)}\to\infty$ for $k<k_{0}$ and ${d}_{n
k}^{(\delta
)}=o_{p}(1)$ for $k\geq k_{0}$. Hence, for any $\theta>0$, as $n\to
\infty$,
\[
P\bigl(\bigl|{d}_{n k}^{(\delta)}\bigr|>\theta\bigr)\to1\qquad \mbox{for }k<k_{0}\quad   \mbox{and}
 \quad
P\bigl(\bigl|{d}_{n k}^{(\delta)}\bigr|>\theta\bigr)\to0\qquad \mbox{for }k\geq k_{0}.
\]

Therefore, for any $\theta>0$ and any $\varepsilon>0$, for each $k$,
there exists a positive integer $N_k$ such that for all $n\geq N_k$,
\[
P\bigl(\bigl|{d}_{n k}^{(\delta)}\bigr|<\theta\bigr)<\varepsilon/(M+1)\qquad \mbox{for
any }
k<k_{0}
\]
and
\[
P\bigl(\bigl|{d}_{n k}^{(\delta)}\bigr|\geq\theta\bigr)<\varepsilon/(M+1)\qquad \mbox{for
any }k_{0}\leq k\leq M.
\]
Note that both $k_0$ and $M$ are finite, we can set an $N$, which is
larger than all~$N_k$ such that the above are satisfied.
Define, for $k\leq M$, \mbox{$B_{n k}:=\{|{d}_{n k}^{(\delta)}|<\theta\}$} and
$B_{n}:= (\bigcap_{i=0}^{k_{0}-1}B_{n i}^{c} )\cup
(\bigcap_{i=k_{0}}^{M}B_{n i} )$ for $n > N$.
Then, for any $\omega\in B_{n}$, $\hat{k}_{\delta,\theta}(\omega)=k_{0}$.
\[
P(B_{n}^{c})\leq\sum_{i=0}^{k_{0}-1}P(B_{n i})+\sum_{k_{0}}^{M}P(B_{n
i}^{c})\leq\varepsilon.
\]
Hence, for any $0<\delta<1$ and $\theta>0$, $\hat{k}_{\delta,\theta}\stackrel{p}{\rightarrow}k_{0}$.

For the case of diverging $k_{0}$, consider $k_{0}<M$ and
$M=o(p^{1/4})$. For any $\theta>0$ and $\delta\leq1/2$, let
$\varepsilon
<\theta/2$ and $\xi>2\theta$. Denote
\begin{eqnarray*}
U_{1}&=& \{\omega\dvtx n^{\delta}|I_{nk}| \leq\varepsilon, \mbox{  for
any $k\leq k_{0}$}  \},\\
U_{2}&=& \{\omega\dvtx n^{\delta}|I_{nk}| \leq\varepsilon, \mbox{  for
any $k_{0}<k\leq M$}  \}
\end{eqnarray*}
 and
 \[
U_{3}= \{\omega\dvtx n^{\delta}(J_{nk}-J_{n k+1})>\xi, \mbox{ for any }
k< k_{0}  \}.\vadjust{\goodbreak}
\]
Then for any $\omega\in\bigcap_{i=1}^{3}U_{i}$, we have $n^{\delta
}(J_{nk}-J_{n k+1})>\xi>2\theta$ for any $k< k_{0}$ and $n^{\delta
}|I_{nk}| \leq\varepsilon<\theta/2$ for any $k\leq M$, which lead to
$n^{\delta}|I_{nk}-I_{n k+1}|<\theta$ for any $k\leq M$. Therefore,
\[
{d}_{n k}^{(\delta)}= n^{\delta}(I_{nk}-I_{n k+1})+n^{\delta
}(J_{nk}-J_{n k+1})>\theta
\qquad\mbox{for any }k< k_{0}
\]
 and
 \[
\big|{d}_{n k}^{(\delta)}\big|\leq n^{\delta}|I_{nk}-I_{n k+1}|<\theta\qquad\mbox{for any }k_{0}\leq k< M.
\]
From (\ref{eq:fixed}), we have $\hat{k}_{\delta,\theta}-k_{0}=0$.
It follows that $\bigcap_{i=1}^{3}U_{i}\subset\{\omega\dvtx \hat
{k}_{\delta
,\theta}-k_{0}=0 \}$. Since $P(\bigcap_{i=1}^{3}U_{i})\to1$ as $n\to
\infty$ by Lemma \ref{lem6}, we have $\hat{k}_{\delta,\theta
}-k_{0}\stackrel{p}{\rightarrow}0$.
\end{pf*}\end{appendix}
\section*{Acknowledgments}
We thank an Associate Editor and two referees
for constructive comments and suggestions which have lead to
improvements in the presentation of the paper. We also thank Professor
Ji Zhu for sharing the data.


%

\printaddresses

\end{document}